\documentclass[11pt]{article}
\newtheorem{Theorem}{Theorem}
\newtheorem{Lemma}{Lemma}
\newtheorem{Corollary}{Corollary}
\newtheorem{Proposition}{Proposition}

\def\qed{\quad \vrule height7.5pt width4.17pt depth0pt}

\newcommand{\lessapprox}{\leq}

\begin{document}

\title{A Singular Parabolic Anderson Model}
\author{Carl Mueller$^1$ and Roger Tribe \\
\\ \\
Dept. of Mathematics \\ University of Rochester \\
Rochester, NY 14627 \\ USA \\
E-mail: cmlr@math.rochester.edu \\ \\
Mathematics Institute \\ University of Warwick \\
Coventry CV4 7AL \\ UK \\
E-mail: tribe@maths.warwick.ac.uk}
\date{}
\maketitle

\begin{abstract}
\noindent
We consider the following stochastic partial differential equation: 
\[
\frac{\partial u}{\partial t}= \frac12 \Delta u+\kappa u \dot F, 
\]
for $x \in \mathbf{R}^d$ in dimension $d \geq 3$, where $\dot F(t,x)$ is a 
mean zero Gaussian noise with the singular covariance 
\[
E\left[\dot F(t,x)\dot F(t,y)\right]=\frac{\delta(t-s)}{|x-y|^2}. 
\]
Solutions $u_t(dx)$ exist as singular measures, under suitable assumptions 
on the initial conditions and for sufficiently small $\kappa$. We 
investigate various properties of the solutions using such tools as 
scaling, self-duality  and moment formulae.  \end{abstract}


\footnotetext[1]{
Supported by an NSF travel grant and an NSA grant.
\par
AMS 1991 \emph{subject classifications} Primary, 60H15; Secondary, 35R60,
35L05.}

\newpage

%
%
\section{Introduction}
%
%
\label{s1} \setcounter{equation}{0}
For readers who want to skip the motivation and definitions, the main results are summarized in Subsection \ref{MainResults}.
\subsection{Background and Motivation}

The parabolic Anderson problem is modeled by the following stochastic partial differential 
equation (SPDE): 
\begin{equation}  \label{1.1}
\frac{\partial u}{\partial t}= \frac12 \Delta u+\kappa u \dot{F}.
\end{equation}
Here $u(t,x) \geq 0$ for $t \geq 0$ and $x \in \mathbf{R}^d$. The equation has various modeling 
interpretations (see Carmona and Molchanov \cite{Carmona+Molchanov94}). The key behavior of 
solutions, called 
intermittency, is that they become concentrated in small regions, often called peaks, separated by 
large almost dead regions.  Except when the covariance of the noise is 
singular at 0,
the linear form of the noise term allows the use of the Feynman-Kac 
formula to study the solutions.  Using this, mostly in the setting of discrete space with a discrete 
Laplacian and with a time-independent noise, there have been many successful descriptions of the
solutions (see \cite{Gartner+Konig+Molchanov00} and the references there to work of Gartner, Molchanov, 
den Hollander, 
Konig and others.) There is less work on the equation with space-time noises but the memoir 
\cite{Carmona+Molchanov94} considers the case of Gaussian noises with various space and time 
covariances.

In addition the ergodic theory of such linear models has been independently studied. Discrete 
versions of the SPDE fit into the framework interacting paricle systems, under the name of 
linear systems.  The reader can consult Liggett \cite{Liggett85}, Chapter IX, Section 4 where, 
using the tools of duality and moments, the ergodic 
behavior of solutions is investigated. This work has been continued for lattice indexed systems 
of stochastic ODEs  (see Cox, Fleischmann and Greven \cite{Cox+Fleischmann+Greven96} and also
Cox,  Klenke and Perkins \cite{Cox+Klenke+Perkins01}). 
The basic picture  is that in dimensions d=1,2 and $d \geq 3$ if $\kappa $ is large, the dead regions get larger 
and  larger and the solutions become locally extinct. Conversely in $d \geq 3$, if $\kappa$ is small, the
diffusion is sufficient to stop the peaks growing and there are non-trivial steady states.

In this paper we study a special case where the noise is white in time and has a space correlation 
that scales, namely 
\begin{equation}  \label{1.2}
E\left[\dot F(t,x)\dot F(s,y)\right]=\frac{\delta(t-s)}{|x-y|^p}.
\end{equation}
The presence of slowly decaying covariances is interesting; one interpretation of the equation given 
in \cite{Carmona+Molchanov94} is in the setting of temperature changes in fluid flow and the noise arises 
as a model 
for the velocities in the fluid, where it is well known that there are slowly decaying covariances (in 
both space and time). Also the equations might arise as a limit of rescaled models where the covariance scaling 
law emerges naturally. Mathematically these covariances are convenient since they imply a scaling 
relation for the solutions that allow us to convert large time behavior into small scale behavior at a
fixed time.

For $0<p<2$ (in dimensions $d \geq 2$) there are function valued solutions with these scaling 
covariances. The Kolmogorov criterion can be used to estimate the H\"{o}lder continuity of solutions 
and in Bentley \cite{Bentley99} the H\"{o}lder continuity is shown to break down as $p \uparrow 2$. In 
this paper we study just the case $p=2$ and establish, in dimensions $d \geq 3$ and when  
$\kappa$ is small, the existence, and uniqueness in law, of measure valued solutions. One can imagine that 
the regularity of solutions breaks down as $p \uparrow 2$ but that there exists a singular, measure 
valued solution at $p=2$ (we do not believe the equation makes sense for the case $p>2$). Note 
that measure valued solutions to an SPDE have been successfully studied in the case of 
Dawson-Watanabe branching diffusions, which can be considered as solutions to the heat equation 
with the noise term $\sqrt{u}dW $, for a space-time white noise $W$ (see Dawson \cite{Dawson93}).

The special covariance $|x-y|^{-2}$ has two singular features: the blow-up near $x=y$ which causes 
the local clustering, so that the solutions become singular measures; and the fat tails at infinity 
which affects large time behavior (for instance we shall prove local extinction in all dimensions). The 
scaling is convenient in that it allows intuition about large time behavior to be  transfered to  results 
on local singularity, and vice-versa. In particular the singularity of the measures can be thought of as a 
description of the intermittency at large  times. 
%
%
\subsection{Definitions}
%
%

Our first task is to give a rigorous meaning to measure valued solutions of  (\ref{1.1}). 
We shall define solutions in terms of a martingale problem.
We do not investigate the possibility of a strong solutions for the equation. We do however 
construct solutions as a chaos expansion with respect to any given noise. These are adapted 
to the same  filtration as the  noise and for some purposes provide a replacement for strong solutions.
One advantage of working with martingale problems is that passing to the limit in approximations 
can be easier with this formulation.

We now fix a suitable state space  for our solutions. Throughout  the paper we consider 
only dimensions $d \geq 3$. The  parameter $\kappa$ will also be fixed to lie in the range 
\begin{equation} \label{1.3}
0 < \kappa <  \frac{d-2}{2}.
\end{equation} 
The restrictions on $d$ and $\kappa$ are due to our requirement that solutions have
finite second moments. We do not explore  the possibility of solutions without
second moments.

Let $\mathcal{M}$ denote the non-negative Radon measures on $\mathbf{R}^d$, 
$\mathcal{C}_c$ the space of continuous functions on $\mathbf{R}^d$ with compact support 
and $\mathcal{C}_c^k$ the space of functions in $\mathcal{C}_c$ with k continuous derivatives. 
We write $\mu(f)$ for the integral $\int f(x) \mu(dx)$ for $\mu \in \mathcal{M}$ and integrable 
$f$, where, unless otherwise indicated, the integral is over the full space $\mathbf{R}^d$. We
consider $\mathcal{M}$ with the vague topology, that is the topology generated by the maps 
$\mu \rightarrow \mu(f)$ for $f \in \mathcal{C}_c$.

The class of allowable initial conditions is described in terms of the singularity of the measures. 
Define 
\[
\|\mu\|_{\alpha}^2 =  \int \! \int \left( 1+ |x-y|^{-\alpha} \right) \mu(dx)\mu(dy)  
\]
and let $\mathcal{H}^a_{\alpha}= \{\mu \in \mathcal{M}: \|\mu(dx) \exp(-a|x|)  \|_{\alpha}<\infty \} $. Note  
the spaces $\mathcal{H}^a_{\alpha}$ are decreasing in $\alpha$ and increasing in $a$. Then define
\[
\mathcal{H}_{\alpha} = \bigcup_{a} \mathcal{H}^a_{\alpha}, \quad 
\mathcal{H}_{\alpha+} = \bigcup_{a} \bigcup_{\beta> \alpha} \mathcal{H}^a_{\beta}, \quad
\mathcal{H}_{\alpha-} = \bigcup_{a} \bigcap_{\beta<\alpha} \mathcal{H}^a_{\beta}.
\]
The sets $\mathcal{H}^a_{\alpha}$ are  Borel subset of $\mathcal{M}$. 
The formula for the second  moments of solutions also leads, for each $d$ and $\kappa$, 
to a distinguished choice of  $\alpha$.
Throughout the paper we make the choice
\[
\alpha = \frac{d-2}{2} -\left[\left(\frac{d-2}{2}\right)^2 - \kappa^2 \right]^{1/2}. 
\]
The restriction (\ref{1.3}) ensures $\alpha\in (0, (d-2)/2)$. We shall require the initial 
conditions to lie in $\mathcal{H}_{\alpha+}$, again to guarantee the existence of second moments.

Suppose $(\Omega, \mathcal{F}, \{\mathcal{F}_t \},P)$ is a filtered probability space. We call an adapted continuous 
$\mathcal{M}$ valued process $\{u_t(dx): t \geq 0 \}$ a (martingale problem) solution to (\ref{1.1}) if it 
satisfies
\begin{enumerate}
\item[i.)] $P(u_0 \in \mathcal{H}_{\alpha+})=1$,
\item[ii.)] $\{ u_t(dx)\}$ satisfies the first and second moment bounds (\ref{1.6}), (\ref{1.7}) given below, and
\item[iii.)] $\{u_t(dx)\}$ satisfies the following martingale problem: for all $f \in \mathcal{C}^2_c$ 
\begin{equation}  \label{1.4}
z_t(f)=u_t(f)-u_0(f) - \int_{0}^{t} \frac12 u_s( \Delta f ) ds
\end{equation}
is a continuous local $\mathcal{F}_t$-martingale with quadratic variation 
\begin{equation}  \label{1.5}
\left\langle z(f) \right\rangle_t = \kappa^2 \int_{0}^{t} \int \! \int  
\frac{f(x) f(y)}{|x-y|^2}u_s(dy)u_s(dx)ds
\end{equation}
\end{enumerate}
If in addition $P(u_0=\mu)=1$ we say that the solution $\{u_t(dx)\}$ has initial condition $\mu$. 

Let $ G_t(x) = (2 \pi t)^{-d/2} \exp( -|x|^2/2t) $.
The moment conditions we require are that  for all $f: \mathbf{R}^d \to [0,\infty)$,
\begin{equation}  \label{1.6}
E \left[ u_t(f)\Big| u_0 \right] =  \int \int  G_t(x-x') f(x') u_0(dx),
\end{equation}
and there exists $C$, depending only on the dimension $d$ and $\kappa$, so that 
\begin{eqnarray}  \label{1.7}
\lefteqn{ E  \left[ \left.\left(  \int f(x) u_t(dx) \right)^2\right|u_0 \right]} \\
& \leq & C  \int_{\mathbf{R}^{4d}} G_{t} (x-x') G_{t}(y-y') f(x') f(y') 
\left(1+  \frac{t^{\alpha} }{|x-y|^{\alpha}|x'-y'|^{\alpha}} \right) 
u_0(dx) u_0(dy) dx' dy'.  \nonumber
\end{eqnarray}
The construction of  solutions in Section \ref{s4} shows that the second moment 
bound is quite natural. We believe that the moment bounds (\ref{1.6}) 
and (\ref{1.7}) are implied by the martingale problem (\ref{1.4}) and (\ref{1.5}), 
although we do not show this. Since establishing second moment bounds is a normal
first step  to finding a solution to the martingale problem, we include these bounds as part of the 
definition of a solution. 

We finish this subsection with some simple consequences of the second moment bound.
\begin{Lemma} \label{momentcor1} 
Suppose $\{u_t(dx)\}$ is a solution to (\ref{1.1}) with initial condition 
$\mu$.  Choose $a$ so that $\mu(dx) \exp(-a|x|) \in \mathcal{H}_{\alpha}$.
\begin{description}
\item[i)] For any $f \in \mathcal{C}_c$ and $t \geq 0$, we have 
\[
E\left[ \int^t_0 \! \int \! \int \frac{f(x) f(y)}{|x-y|^2} u_s(dx) u_s(dy) ds \right] < \infty 
\]
and hence the process $z_t(f)$ defined in (\ref{1.4}) is a true martingale. 
\item[ii)] For any $0 \leq \rho < d-\alpha$ and  $t >0$ 
\[
E\left[ \int \! \int \left( 1 + |x-y|^{-\rho}\right) e^{-a|x|-a|y|} u_t(dx) u_t(dy) \right] < \infty 
\]
and hence $u_t \in \mathcal{H}_{(d- \alpha)-}$ almost surely.
\end{description}
\end{Lemma}
\noindent
\textbf{Proof.} For part i) it is sufficient to check that $E[ \langle z(f) \rangle_t] < \infty$ to ensure 
that $z_t(f)$ is a true martingale. Using the second moments (\ref{1.7}) we have 
\begin{eqnarray}
\lefteqn{ E\left[ \int^t_0 \int \! \int \frac{f(x) f(y)}{|x-y|^2} u_s(dx) u_s(dy) ds \right] }  \label{1.8} \\
& \leq  & C  \int^t_0 \int_{\mathbf{R}^{4d}} G_s(x-x') G_s(y-y')  \frac{f(x')f(y')}{|x'-y'|^2} 
 \left(1+ \frac{s^{\alpha}}{ |x-y|^{\alpha} |x'-y'|^{\alpha}} \right) 
\mu (dx) \mu (dy) dx' dy' ds. \nonumber
\end{eqnarray}
We now estimate the $dx' dy'$ integral in the above expression by using the simple bound, for $0 \leq r < d$, 
\begin{equation}  \label{1.9}
\int \int G_t(x-x') G_t(y-y') |x'-y'|^{-r}dx' \, dy' \leq C(r ) \left( |x-y|^{-r} \wedge t^{-r/2} \right).
\end{equation}
For compact support $f$ and any $a>0$ we have the bound 
\begin{equation} \label{1.10}
\int G_s(x-x') f(x') dx' \leq C(a,f,t) e^{-a|x|} \quad \mbox{for all $s \leq t$, $x \in \mathbf{R}^d$}.
\end{equation}
Then, applying H\"{o}lder's inequality with $1<p<d/2$ and $p^{-1} + q^{-1} =1$,
we have, for all $s \leq t$, 
\begin{eqnarray*}
\lefteqn{ \int \int G_s(x-x') G_s(y-y') \frac{f(x') f(y')}{|x'-y'|^2} dx' dy' } \\
& \leq &  \left(\int \int G_s(x-x') G_s(y-y') \frac{1}{|x'-y'|^{2p}} dx' dy' \right)^{1/p} \\
&& \hspace{.5in} \cdot  \left(\int \int G_s(x-x') G_s(y-y') f^q (x') f^q (y') dx' dy' \right)^{1/q} \\
& \leq & C(a,f,t) e^{-a|x| - a|y|} \;  ( |x-y|^{-2} \wedge s^{-1}) \\
& \leq & C(a,f,t) e^{-a|x| - a|y|} \; |x-y|^{-\alpha \wedge 2} \; s^{- (2-\alpha)_+/2} \\
& \leq & C(a,f,t) e^{-a|x| - a|y|} \; \left( 1+ |x-y|^{-\alpha} \right) \; s^{-(2-\alpha)_+/2}. 
\end{eqnarray*}
A similar calculation, using $2+\alpha < d$, gives the bound 
\[
\int \int G_s(x-x') G_s(y-y')  \frac{f(x') f(y')}{|x'-y'|^{2+\alpha}} dx' dy'
 \leq  C(a,f,t)  s^{-(2+\alpha)/2}. 
\]
Now we substitute these bounds into (\ref{1.8})  to obtain 
\begin{eqnarray*}
\lefteqn{E\left[ \int \int \frac{f(x) f(y)}{|x-y|^2} u_s(dx) u_s(dy) ds \right]} \\
&  \leq & C(a,f,t) \int^t_0 s^{-(2-\alpha)_+/2}  \int \int \left(1+ |x-y|^{-\alpha}\right)
e^{-a|x| -a|y|} \mu (dx) \mu (dy)  ds.
\end{eqnarray*}
which is finite since $\|\mu(dx) \exp(-a|x|) \|_{\alpha}< \infty$.  

For part ii) use the second moment bound (\ref{1.7}) to see that
\begin{eqnarray}
\lefteqn{ E \left[ \int \! \int \left( 1+ |x-y|^{-\rho} \right) 
e^{-a|x|-a|y|} u_t(dx) u_t(dy) \right] }  \label{1.11}  \\
& \leq & C \int_{\mathbf{R}^{4d}} \! G_t(x-x')G_t(y-y') 
\left( 1+ |x'-y'|^{-\rho} \right) e^{-a|x'|-a|y'|} \nonumber \\
&& \hspace{.5in} \cdot  
\left(1+ \frac{t^{\alpha}}{|x-y|^{\alpha}|x'-y'|^{\alpha}}\right) 
\mu(dx) \mu(dy) dx' dy'. \nonumber  
\end{eqnarray}
Using the bound 
$\int G_t(x-x') \exp(-a|x'|) dx' \leq C(t,a) \exp(-a|x|)$ and (\ref{1.9}),
 we estimate the $dx'dy'$ integral in a similar manner as above. 
We illustrate this only on the most singular term. 
For  $p,q>1$ with $p^{-1}+q^{-1}=1$, 
\begin{eqnarray}
\lefteqn{ \int \! \int G_t(x-x') G_t(y-y') |x'-y'|^{-(\rho +\alpha)} 
e^{-a|x'|-a|y'|} dx' dy' } \nonumber \\
& \leq &  \left(\int \int G_t(x-x') G_t(y-y') e^{-ap|x'|-ap|y'|} dx' dy' \right)^{1/p} \nonumber \\
& & \hspace{.5in} \cdot 
\left(\int \!\int G_t(x-x') G_t(y-y') |x'-y'|^{-q(\rho + \alpha)}dx' dy' \right)^{1/q} \nonumber \\
& \leq & C(a, \rho,p,q) e^{-a|x| - a|y|} \left( |x-y|^{-(\rho+\alpha)} \wedge t^{-(\rho+\alpha)/2} \right) 
\label{1.12} \\
& \leq & C(t,a,\rho,p,q) e^{-a|x| - a|y|}. \nonumber  
\end{eqnarray}
provided that $q(\rho + \alpha)< d$. Such a $q>1$ can be found whenever 
$\rho+\alpha <d$. Substituting this estimate into (\ref{1.11}) gives the result.  \qed
%
%
%
\subsection{Main Results} \label{MainResults}
%
%

We start with a result on existence and uniqueness. 
\begin{Theorem} \label{t1} 
For any $\mu \in \mathcal{H}_{\alpha+} $ there exists a solution to (\ref{1.1}) started 
at $\mu$. Solutions starting at $\mu \in \mathcal{H}_{\alpha+}$
are unique in law. If we denote this law by   $Q_{\mu}$ then 
the  set $\{Q_{\mu}: \mu \in \mathcal{H}_{\alpha+}\}$ forms a Markov family of laws. 
\end{Theorem}
The existence part of Theorem \ref{t1} is proved in Section \ref{s4} and the uniqueness in Section 
\ref{s5}. 
The  next theorem, which is proved in Section \ref{s6}, shows death from finite initial conditions 
and  local extinction from certain infinite initial conditions. Write $B(x,r)$ for the open ball or radius $r$ 
centered at $x$. We say that a random  measure $u_0$ has bounded local intensity if 
$E[ u_0(B(x,1))]$ is a bounded function of $x$.
\begin{Theorem} \label{t4} 
Suppose $\{u_t(dx)\}$ is a solution to (\ref{1.1}).
\begin{description}
\item[i) Death from finite initial conditions.]  If  $P(u_0 \in \mathcal{H}_{\alpha}^0)=1 $  then 
$(u_t,1) \rightarrow 0 $ almost surely as $t \rightarrow \infty$.
\item[ii) Local extinction from infinite initial conditions.]  If $u_0$ has bounded local intensity 
and $A \subseteq \mathbf{R}^d$ is a bounded set then $u_t(A) \rightarrow 0$ in probability 
as $t \rightarrow \infty$.
\end{description}
\end{Theorem}

Finally, we state our main results describing the nature of the measures $u_t(dx)$. These are proved in 
Section \ref{s7}.
\begin{Theorem} \label{t5} 
Suppose that $\{u_t(dx)\}$ is a solution to (\ref{1.1}) satisfying 
$P( u_0 \neq 0 )=1$.  Fix $t>0$.  Then the following properties 
hold with probability one.
\begin{description}
\item[i) Dimension of support.]   If a Borel set $A$ supports the measure $u_t(dx)$  then the 
Hausdorff dimension of $A$ is at least $d-\alpha$.
\item[ii) Density of support.]  The closed support of $u_t(dx)$ is $\mathbf{R}^d$.
\item[iii) Singularity of solutions.]  The absolutely continuous part of $u_t(dx)$ is zero.
\end{description}
\end{Theorem}
\textbf{Remarks}

\textbf{1.}  Although Theorem \ref{t5} gives an almost sure result for fixed $t$, it leaves open 
the possibility that  there are random times at which the properties fail. In Section \ref{s7} we 
shall show that $P(u_t \in \mathcal{H}_{\alpha+} \; \mbox{for all $t \geq 0$} )=1$. 
This  implies that the  weaker lower bound $d-2-\alpha$ on the dimension of supporting sets 
is valid for all times. 

\textbf{2.} The reader might compare 
the behavior described in Theorem \ref{t5}  with that of the Dawson-Watanabe branching diffusion 
in $\mathbf{R}^d$, for $d \geq 2$. This is a singular measure valued process whose support is two
dimensional,  and, if started with a finite measure of compact support,  has compact support for all time. 

\textbf{3.} Many of the results go through for the boundary case $\kappa = (d-2)/2$ and
for initial conditions in $\mathcal{H}_{\alpha}$, although we have not stated results in these cases.
The chaos expansion in Section \ref{s4} holds in both these boundary cases and the second moments are finite. Although our proof that the chaos expansion
satisfies (\ref{1.5}) uses $\kappa <(d-2)/2$ and $\mu \in \mathcal{H}_{\alpha+}$ we do not believe these restrictions are needed for this. However our proof of uniqueness for solutions in Section \ref{s5} does seem to require the strict inequalities. This leaves open the possibility that there are solutions with a different law to that constructed via the chaos expansion. 
Theorems \ref{t4} and \ref{t5} will hold in the boundary cases for the solutions constructed via chaos expansion solutions. Parts of Theorems \ref{t4} and \ref{t5} also hold for all solutions, for example Propositions \ref{liggett} and \ref{dimension} use only the martingale problem in their proof and  hold for any solution in the boundary cases.
%
%
\subsection{Tools} \label{Tools}
%
%

We briefly introduce the main tools that we use. The first tool, simple scaling for the equation, is summarized  
in the following lemma.
\begin{Lemma}
\label{scaling} Suppose that $\{u_t(dx)\}$ is a solution to (\ref{1.1}). Let  
$a,b,c>0$ and define 
\[
v_t(A)=au_{bt}(cA) \quad \mbox{for Borel $A \subseteq \mathbf{R}^d$}
\]
where $cA= \{cx:x \in A\}$. Then $\{v_t(dx)\}$ is a solution to the equation 
\[
\frac{\partial v_t}{\partial t} = \frac{b}{2c^2}\Delta v + \kappa 
\frac{b^{1/2}}{c} v \dot{F}_{b,c}(t,x) 
\]
where $\dot F_{b,c}(t,x)$ is a Gaussian noise identical in law to $\dot{F}(t,x)$.
\end{Lemma}
The equation for $\{v_t(dx)\}$ is interpreted via a martingale problem, as in (\ref{1.1}).
The easy proof of this lemma is omitted.
 
The next tool is our equation for the second moments.  The linear noise term implies 
that the solutions have closed moment equations. By this we mean that the moment densities 
\[
H_t(x_1,x_2, \ldots,x_n)dx_1\dots dx_n = E\left[u_t(dx_1) u_t(dx_2)\ldots u_t(dx_n)\right] 
\]
satisfy an autonomous PDE. Formally assuming the solution has a smooth density $u_t(x)$, 
applying Ito's formula to the product $u_t(x_1)\dots u_t(x_n)$ and taking expectations suggests 
that $H_t$ satisfies 
\[
\frac{\partial H_t}{\partial t} = \frac12  \Delta H_t + \kappa^2 H_t  
\sum_{1 \leq i < j \leq n} \frac{1}{|x_i - x_j |^2} . 
\]
Then the Feynman-Kac representation for this linear equation suggests that 
\[
H_t(x_1, \ldots, x_n) =  E_{x_1,\ldots,x_n} \left[ u_0 ( X^1_t) \ldots u_0( X^n_t) \exp 
\left(\int^t_0 \sum_{1 \leq i < j \leq n} \frac{\kappa^2}{ |X_s^i - X^j_s|^{2}} ds \right) \right] 
\]
where $ E_{x_1,\ldots,x_n} $ denotes expectation with respect to $ n $ independent d-dimensional Brownian 
motions. This formula makes sense when $u_0$ has a density, but more generally we can
expect for solutions $\{u_t(dx)\}$ to (\ref{1.1}) started at $\mu$, and when 
$f_i \in \mathcal{C}_c$ for $i=1,\ldots,k$, 
\begin{eqnarray}
\lefteqn{ E\left[\prod_{i=1}^{n}u_t(f_i)\right] } \label{1.13} \\
& = &  \int_{\mathbf{R}^{2nd}}  E_{0,x_1,\dots,x_n}^{t,y_1,\dots,y_n}
\left[\exp\left(\sum_{1\le j<k\le n}\int_{0}^{t} \frac{\kappa^2} {\left|X^{(j)}_s-X^{(k)}_s \right|^{2}}
 ds \right) \right]  \prod_{i=1}^{n}G_t (x_i-y_i) f_i (y_i) \mu (dx_i) dy_i  \nonumber
\end{eqnarray}
where $ E_{0,x_1,\dots,x_n}^{t, y_1,\dots,y_n} $
is expectation with respect to $ n $ independent d-dimensional Brownian bridges 
$ (X^{(1)}_t, \ldots, X^{(n)}_t) $ started at $(x_i)$ at time zero and ending at 
$ (y_i)$ at time $ t $. In Section \ref{s2} we investigate the values of $ \kappa $ for which this 
expectation is finite.

The next tool is the expansion of the solution a Wiener chaos expansion, involving multiple integrals 
over the noise $F(t,x)$. Wiener chaos expansions have been used before for linear 
equations; for example see Dawson and Salehi \cite{Dawson+Salehi80} or Nualart and Zakai 
\cite{Nualart+Zakai89}. The 
idea is to start with the Green's function representation, assuming (falsely) that a function 
valued solution exists:
\begin{equation}  \label{1.14}
u_t (y) = G_t \mu (y)  + \kappa \int^t_0 \int G_{t-s} (y-z) u_s (z) F(dz,ds).
\end{equation}
The first term on the right hand side of this representation uses the notation 
$ G_t \mu (y)=\int G_t(y-z) \mu(dz)$. The second term involves again the the non-existent density 
$u_s(z)$. However we can use the formula for $u_t(y)$ given in (\ref{1.14}) to substitute for the term 
$u_s(z)$ which appears on its right hand side. The reader can check that if we keep repeating this 
substitution, and assume the remainder term vanishes, we will arrive at the following formula: 
for a test function $f \in \mathcal{C}_c$, 
\begin{equation}  \label{1.15}
u_t (f) = \sum_{n=0}^{\infty}I^{(n)}_t (f, \mu)
\end{equation}
where
\begin{equation} \label{1.16}
 I^{(n)}_t(f,\mu) = \int \! \int f(y) I^{(n)}_t(y,z) \mu(dz) dy
\end{equation}
and where the $I^{(n)}$ are defined as follows: $I^{(0)}_t(y,z) = G_t(y-z)$ and for $n \geq 1$ 
\begin{equation}  \label{1.17}
I^{(n)}_{s_{n+1}} (y_{n+1}, z) = \kappa^n \int_0^{s_{n+1}} \int^{s_n}_0 {\dots} \int_{0}^{s_2} 
\int_{\mathbf{R}^{nd}} G_{s_1}(y_1-z) \prod_{i=1}^{n}G_{s_{i+1}-s_i}(y_{i+1}-y_i) F(dy_i,ds_i).
\end{equation}
In Section \ref{s4} we shall show  that the stochastic integrals in (\ref{1.17}) are well defined, and the 
series $(\ref{1.15})$  converges in $L^2$ and defines a solution. The point is that the series
$\sum_n I^{(n)}_t(y,z)$ does not converge pointwise, but after smoothing by integrating against the initial
measure and the test function the series does converge. The restriction (\ref{1.3}) 
on $\kappa$ and the choice of space $\mathcal{H}_{\alpha+}$ for the initial conditions is exactly
what we need to ensure this $L^2$ convergence. For larger values of $\kappa$ it is possible
that the series converges in $L^p$ for some $p<2$. It is also always possible to consider the
chaos expansion (\ref{1.15}) itself as a solution, if we interpret solutions in a suitably
weak fashion, for example as a linear functionals on Wiener space. We do not investigate either 
of these possibilities. 

The symmetry of the functions $I^{(n)}_t(y,z)$ in $y$ and $z$ makes it clear that a time reversal
property should hold. This is well known for linear systems and for the parabolic Anderson model, 
and is often called self duality. Suppose that $\{u_t(x)\}, \, \{v_t(x)\}$ are two solutions of (\ref{1.1}) 
started from suitable absolutely continuous initial conditions $u_0(x)dx$ and $v_0(x)dx$. We
expect that $ u_t (v_0) $ has the same distribution as $ u_0 (v_t) $. In Section \ref{s5} we shall 
use this equality to establish uniqueness of solutions.

The Feynman-Kac formula is a standard tool in analogous discrete space models. In the continuous 
space setting of the parabolic Anderson equation (\ref{1.1}), we shall replace the noise $F$ by a 
noise $\bar{F}$ that is Gaussian, white in time and with a smooth, translation invariant covariance 
$\Gamma (x-y)$ in space. Then the Feynman-Kac representation is 
\begin{eqnarray}  
u_t(x) & = & E_x \left[u_0(X_t ) e^{-\Gamma(0)t} \exp \left( \kappa \int_{0}^{t}
 \bar{F}(ds,X_{t-s}) \right) \right]  \nonumber \\
& = & \int G_{t} (x-y) u_0(y) e^{-\Gamma (0) t}  E_{0,y}^{t,x} \left[ 
\exp \left( \kappa \int_{0}^{t} \bar{F} (ds, X_{s}) \right) \right]. \label{1.18}
\end{eqnarray}
A proof of this representation can be found in Kunita \cite{Kunita90} 
Theorem 6.2.5 and we make use of it in Section \ref{s7}. Since our covariance blows 
up at the origin the exponential factor $\Gamma(0)$ is infinite 
and the representation can only be used for approximations.

Finally a remark on notation: we use $C(t, p, \ldots) $ for a constant whose exact value is unimportant 
and may change from line to line, that may depend on the dimension $d$ and the parameter $\kappa $ 
(and hence also on $\alpha$),  but whose dependence on other parameters will be indicated. 
%
%
\section{A Brownian exponential moment}
%
%
\label{s2}  \setcounter{equation}{0} 

As indicated in the introduction, the second moments of solutions $\{u_t(dx)\}$ to (\ref{1.1}) can be 
expressed in terms of the expectation of a functional of a Brownian bridge. An upper bound for 
these expectations is a key estimate in the construction of our solutions. In this section we show the 
following bound.
\begin{Lemma} \label{exponentialmoment}
For all $0 \leq \eta \leq (d-2)^2/8$ there exists $C(\eta) < \infty$ so that for all $x,y,t$
\[
 E^{t,y}_{0,x} \left[ \exp \left( \eta \int^t_0 \frac{ds}{|X_s|^2} \right) \right]  \leq 
C(\eta) \left(1+ \frac{t}{|x|\, |y|}\right)^{\alpha(\eta)}
\]
where 
\[ 
\alpha(\eta) = \frac{d-2}{2}  - \left[\left(\frac{d-2}{2} \right)^2 - 2 \eta \right]^{1/2}. 
\]
\end{Lemma}
We first treat the case of Bessel processes and Bessel bridges (see Revuz and Yor \cite{Revuz+Yor91} 
chapter XI for the basic definitions). The reason for this is that the laws of two Bessel processes, of two 
suitable different dimensions, are mutually absolutely continuous and the Radon-Nikodym derivative 
involves exactly the exponential functional we wish to estimate. 

Let $C[0,t]$ be the space of real valued continuous paths up to time $t$ and let $\{R_t\}$ be the canonical 
path variables. For $d \in [2,\infty)$ and $a,b > 0$ we write $E^{(d)}_a$ for expectations under the law of  the 
$d$-dimensional Bessel process started at $a$ and   $q^{(d)}_t(a,b)$ for the transition density. We write 
$E^{(d)}_{a,b,t}$ for expectations under the law of the $d$-dimensional  Bessel bridge starting at $a$ and 
ending at $b$ at time $t$. Suppose that $Y$ is a non-negative random variable on the space $C[0,t]$, 
measurable with respect to $\sigma (R_s: s \leq t)$. Lemma 4.5 of Yor \cite{Yor80}, (or Revuz and Yor 
\cite{Revuz+Yor91}, Chapter XI, exercise 1.22), expressed in our notation, states that the following 
relationship holds: if 
$\lambda, \mu \geq 0$ then 
\begin{equation}  \label{2.1}
E_a^{(2 \lambda+2)} \left[Y \exp \left(-\frac{\mu^2}{2} \int_{0}^{t} \frac{ds}{R_s^2} \right)
\left(\frac{R_t}{a}\right)^{-\lambda}\right]  = E_a^{(2\mu+2)}\left[Y\exp \left(-\frac{\lambda^2}{2}
\int_{0}^{t} \frac{ds}{R_s^2} \right) \left(\frac{R_t}{a}\right)^{-\mu}\right].
\end{equation}
Now for $0 \leq \eta \leq (d-2)^2/8$ we 
choose values for $\lambda, \mu, Y$ in this identity as follows: 
\[ 
\lambda = \frac{d-2}{2}, \quad \mu = \left[\left(\frac{d-2}{2} \right)^2- 2 \eta \right]^{1/2},
\quad Y = \exp \left( \left[ \eta + \frac{\mu^2}{2}\right] \int_{0}^{t} 
\frac{ds}{R_s^2} \right) \left( \frac{R_t}{a} \right)^{\lambda}
\mathbf{1}(R_t \in db).
\]
Note with these choices that $\alpha(\eta) = \lambda - \mu $, $ 2 \eta + \mu^2 - \lambda^2 = 0 $, and
$d=2\lambda+2$. Applying (\ref{2.1}) we find 
\begin{eqnarray*}   
\lefteqn{  E_a^{(d)} \left[ \exp \left( \eta \int_{0}^{t} \frac{ds}{R_s^2} \right) 
\mathbf{1} (R_t \in db) \right] } \\
&=& E_a^{(2 \lambda + 2)}\left[ Y \exp \left( - \frac{\mu^2}{2} \int_{0}^{t} 
\frac{ds}{R_s^2} \right)  \left( \frac{R_t}{a} \right)^{-\lambda}\right] \\
&=& E_a^{(2\mu+2)} \left[Y \exp \left( - \frac{\lambda^2}{2} \int_{0}^{t} 
\frac{ds}{R_s^2} \right) \left( \frac{R_t}{a} \right)^{-\mu} \right] \\
&=& E_a^{(2\mu+2)} \left[ \exp \left( \left[ \eta + \frac{\mu^2}{2} -
\frac{\lambda^2}{2} \right]  \int_{0}^{t} \frac{ds}{R_s^2} \right)
 \left(\frac{R_t}{a} \right)^{\lambda - \mu}
\mathbf{1} (R_t \in db)   \right] \\
&=& E_a^{(2\mu+2)} \left[ \left( \frac{R_t}{a} \right)^{\alpha(\eta)}
\mathbf{1} (R_t \in db)  \right] \\
&=& a^{- \alpha(\eta)} b^{\alpha(\eta)} q^{(2\mu+2)}_t(a,b) db. 
\end{eqnarray*}
Hence
\begin{equation} \label{2.2}
 E_{a,b,t}^{(d)} \left[ \exp \left( \eta \int_{0}^{t} \frac{ds}{R_s^2} \right)  \right]
= a^{- \alpha(\eta)} b^{\alpha(\eta)} 
\frac{q^{(2\mu+2)}_t(a,b)}{q^{(d)}_{t}(a,b)}. 
\end{equation}
There is an exact formula for the Bessel transition density
\[
 q^{(d)}_t(a,b) = t^{-1} a^{-(d-2)/2}b^{d/2} \exp(-(a^2+b^2)/2t) I_{(d/2)-1}(ab/t) 
\]
in terms of the (modified) Bessel functions $I_{\nu}$ of index $\nu = (d/2)-1$. The  Bessel functions
$I_{\nu}(z)$ are continuous and strictly positive for $z \in (0,\infty)$ and satisfy the asymptotics,
for $c_1, c_2>0$,
\[
I_{\nu}(z) \sim c_1 z^{\nu} \quad \mbox{ as $ z \downarrow 0$,}  \quad I_{\nu}(z) \sim c_2 z^{-1/2} 
e^{z} \quad \mbox{ as $ z \uparrow \infty$.}  
\]
Using these we find that
\begin{equation} \label{2.3}
E_{a,b,t}^{(d)} \left[ \exp \left( \eta \int_{0}^{t} \frac{ds}{R_s^2} \right)  \right]
\leq C( \eta) \left( 1+ \frac{t}{ab} \right)^{\alpha(\eta)} \quad
\mbox{for all $a,b,t>0$.}
\end{equation}

We now wish to obtain a similar estimate for a Brownian bridge. Recall the skew product 
representation for a $d$-dimensional Brownian motion $X_t$, started from $x \neq 0$. There 
is a Brownian motion $W(t)$ on the sphere $\mathbf{S}^{d-1}$, started at $x/|x|$ and independent
of $X$, so that 
\[ 
X_t/|X_t| = W\left( \int^t_0 |X_s|^{-2} ds\right).
\]
We may find a constant $C$ so that
$  P_x( W(t) \in d\theta) \leq C d \theta $ for all $x \in \mathbf{S}^{d-1}$ and $t \geq 1$. 
We now consider the exponential moment for a $d$-dimensional  Brownian bridge running 
from $x \neq 0$ to $y \neq 0$ in time $t$.
\begin{equation} \label{2.4}
E^{y,t}_{x,0} \left[ \exp \left( \eta \int^t_0  \frac{ds}{|X_s|^2} \right) \right] 
 \leq e^{\eta} + E^{y,t}_{x,0} \left[ \exp \left( \eta \int^t_0  \frac{ds}{|X_s|^2} \right)  
\mathbf{1}  \left( \int^t_0 \frac{ds}{|X_s|^2} \geq 1 \right) \right].
\end{equation}
Now we estimate the second term on the right hand side of (\ref{2.4}). 
\begin{eqnarray*}
\lefteqn{ E^{y,t}_{x,0} \left[ \exp \left( \eta \int^t_0 
\frac{ds}{|X_s|^2} \right)  \mathbf{1} 
\left( \int^t_0 \frac{ds}{|X_s|^2} \geq 1 \right) \right]  } \\
& = & \frac{1}{G_t(x-y)} E_x  \left[ \exp \left( \eta \int^t_0 
\frac{ds}{|X_s|^2} \right)  \mathbf{1} 
\left( \int^t_0 \frac{ds}{|X_s|^2} \geq 1, X_t \in dy \right) \right]  \\
& = & \frac{C |y|^{1-d}}{G_t(x-y)}
E_x \left[ \exp \left( \eta \int^t_0 
\frac{ds}{|X_s|^2} \right)  \mathbf{1} \left(  \int^t_0 
\frac{ds}{|X_s|^2} \geq 1, |X_t| \in d|y|, W\left(\int^t_0 
\frac{ds}{|X_s|^2}\right) \in d(y/|y|) \right) \right] \\
& \leq &\frac{C |y|^{1-d}}{G_t(x-y)} 
E_x \left[ \exp \left( \eta \int^t_0 
\frac{ds}{|X_s|^2} \right)  \mathbf{1} \left(  \int^t_0 
\frac{ds}{|X_s|^2} \geq 1, |X_t| \in d|y|  \right) \right] \\
& \leq &\frac{C |y|^{1-d} q^{(d)}_t(|x|,|y|)}{G_t(x-y)} 
E^{(d)}_{|x|,|y|,t}  \left[ \exp \left( \eta \int^t_0 
\frac{ds}{|R_s|^2} \right) \right].
\end{eqnarray*}
Using the explicit representation for the Bessel density given above we find that
\[
\frac{|y|^{1-d}q^{(d)}_t(|x|,|y|)}{G_t(x-y)}  \leq C(R ), \quad \mbox{whenever }\frac{|x| |y|}{t} \leq R.
\] 
Combining this with  (\ref{2.4}) and our estimate (\ref{2.3}) for the 
Bessel bridge we obtain the desired bound for $(x,y,t)$ in any region where
 $ \{ |x| |y|/t  \leq R\} $. 

We felt there should be a short way to treat the remaining case, but we seem to need a slightly complicated argument to treat the case $|x| |y|/t $ large. Note our aim is only to find a constant bound for the exponential moment in this region. Define
\[
F_K(x,y,t) := E^{t,y}_{0,x} \left[ \exp \left( \eta \int^t_0 \frac{ds}{|X_s|^2} \wedge K \right) \right]
\leq E^{t,y}_{0,x} \left[ \exp \left( \eta \int^t_0 \frac{ds}{|X_s|^2} \right) \right]
=: F(x,y,t).
\]
Brownian scaling implies that $F(x,y,t) = F(c^{1/2}x, c^{1/2}y,ct)$ for any $c>0$. So we 
may scale time away and it is enough to control $F(x,y,1) $. 
We have proved above, for any $R$,  
\begin{equation} \label{2.5}
F(x,y,1)  \leq  C(R, \eta) \left(1+ \frac{1}{|x|\, |y|}\right)^{\alpha(\eta)}  
 \mbox{whenever $|x| |y| \leq R$.}
\end{equation}
We first show we may reduce to the case where $|x|=|y|$. Suppose that $|x||y| \geq 1$ and  $|x|>|y|$. 
Define stopping times
\[ 
\sigma_1 = \inf  \{ t: |X_t| \leq |y| \}, \quad \sigma_2 = \inf  \{ t: |X_t| |y| \leq 1-t \},
\]
and let $ \sigma= \sigma_1 \wedge \sigma_2$.
Note that for $t < \sigma_1$ we have $1/|X_t| \leq 1/|y|$ and for 
$t < \sigma_2$ we have $1/|X_t| \leq |y|/(1-t)$. So we can bound the 
integral in $F(x,y,1)$ by

\newpage

\begin{eqnarray*}
\int^1_0 \frac{dt}{|X_t|^2} & \leq & \int^{\sigma}_0 \frac{dt}{|X_t|^2}
+ \int^1_{\sigma} \frac{dt}{|X_t|^2} \\
& \leq & \int^{\sigma}_0 \left( \frac{|y|^2}{(1-t)^2} \wedge \frac{1}{|y|^2}  \right) dt
+ \int^1_{\sigma} \frac{dt}{|X_t|^2} \\
& \leq & \int^{(1-|y|^2)_+}_0  \frac{|y|^2}{(1-t)^2} dt
+ \int^1_{(1-|y|^2)_+} \frac{1}{|y|^2} dt
+ \int^1_{\sigma} \frac{dt}{|X_t|^2} \\
& \leq & 2 + \int^1_{\sigma} \frac{dt}{|X_t|^2}.
\end{eqnarray*}
Conditioned on the values of $\sigma$ and $X_{\sigma}$, the path between $t \in [\sigma,1]$
is a new Brownian bridge. Hence
\[
F(x,y,1)   \leq  e^{2 \eta}   E \left(  F(X_{\sigma},y, 1 -\sigma) \right)  
= E \left[ F \left(\frac{X_{\sigma_1}}{(1-\sigma_1)^{1/2}}, \frac{y}{(1-\sigma_1)^{1/2}}, 1\right) \right].
\]
By definition $|X_{\sigma_2}| | y| / ( 1- \sigma_2) = 1 $ so that
$ F (X_{\sigma_2}/ (1-\sigma_2)^{1/2}, y/ (1-\sigma_1)^{1/2}, 1)$ can be  bounded by a 
constant using (\ref{2.5}). Also  
on the set $\{\sigma_1 < \sigma_2 \}$ we know that $|X_{\sigma_1}| | y| / ( 1- \sigma_1) \geq 1 $
and $|X_{\sigma_1}| = |y|$. So 
if we can bound the $F(x,y,1)$ on the set diagonal case $\{ |x|=|y|, |x| |y| \geq 1 \} $ we can bound
$ F (X_{\sigma_1}/ (1-\sigma_1)^{1/2}, y/ (1-\sigma_1)^{1/2}, 1)$,  and in consequence also $F(x,y,1)$.

We now give a brief sketch to motivate the final argument. 
Consider the ``worst case'' of a bridge from $x=Ne_1$ to  $y=-Ne_1$ over  time one. Run both ends of the bridge until both ends first hit  the ball of radius $N/2$. 
When $N$ is large the bridge will enter the ball near $x/2$ and exit near $y/2$ and spend 
close to time $1/2$ inside the ball.  We may therefore  approximately bound the exponential as
\[ 
\exp \left( \int^1_0 \frac{1}{|X_s|^2} ds \right)  \lessapprox \exp( 4 N^{-2}) \exp \left( \int^{1/4}_{3/4} \frac{1}{|X_s|^2} ds \right).
\]
Using the scaling of $F(x,y,t)$ we see that $F(Ne_1, -Ne_1,1)$ is approximately bounded by 
\[
\exp(4 N^{-2}) F(Ne_1 /2, -N e_1 /2, 1/2) =  \exp(4N^{-2}) F(Ne_1 /2^{1/2}, -N e_1 /2^{1/2}, 1).
\] 
By iterating this argument we will bound $F(Ne_1,-Ne_1,1)$ for large $N$ by values for 
small $N$ where we know it is bounded by (\ref{2.5}). 

We now give the basic iterative construction. Suppose 
that $|x|=|y|=R \geq 1$ and consider the Brownian bridge
$\{X_t\}$ from $x$ to $y$ in time $1$. Define random times
\[
\sigma= \inf \{t: |X_t| \leq R/2\}, \quad \tau= \sup \{t: |X_t| \leq R/2\}
\]
on the set $\{ \inf_t |X_t| < R/2\}= \{ \sigma < \tau\}$. On the 
set $\{ |X_t| > R/2, \, \forall t \in [0,1]\}$ we have the bound
$ \int^1_0 |X_s|^{-2} ds \leq 4 R^{-2}$. On $\{ \sigma < \tau\} $ we have the bound
\[
\int^1_0 \frac{dt}{|X_t|^2}  \leq 4 R^{-2}  \int^{\tau}_{\sigma} \frac{dt}{|X_t|^2}.
\]
Conditioned on $\sigma, \tau, X_{\sigma}, X_{\tau}$, the part of the path
$\{X_t: t \in [\sigma,\tau]\}$ is a new Brownian bridge. So we may estimate
\begin{equation}
F(x,y,1) \leq  \exp( 4 \eta R^{-2}) \left( P(\{ |X_t| > R/2, \, \forall t \in [0,1]\} ) + 
 E \left( F(X_{\sigma}, X_{\tau}, \tau-\sigma) \mathbf{1} (\sigma< \tau)  \right) \right).
\label{2.6}
\end{equation}
The same bound holds with $F$ replaced by $F_K$. 

We will repeat this construction with a new Brownian bridge running from 
$ X_{\sigma}/ (\tau- \sigma)^{1/2}$ to $ X_{\tau}/ (\tau- \sigma)^{1/2}$. The following lemma
shows that when $R$ is large we have usually made an improvement in that
this bridge is closer to the origin.
\begin{Lemma} \label{newlemma} 
There exists $\gamma <1$ and $c_3< \infty, c_4>0$ so that, when $|x|=|y|=R$, 
\begin{equation} \label{2.7}
P\left( \sigma < \tau, \; X_{\sigma} \cdot X_{\tau} \leq 0, \;
 \left| X_{\sigma}/(\tau- \sigma)^{1/2} \right| \geq \gamma |x| \right)
\leq c_3 \exp(-c_4 R),
\end{equation}
and there exist $c_5< \infty, c_6 >0$ so that, if in addition $x \cdot y \geq 0$,
\begin{equation} \label{2.8}
P\left(  \inf_t |X_t| < R/2 \right)
\leq c_5 \exp(-c_6 R),
\end{equation}
\end{Lemma}

\noindent
\textbf{Proof}. 
We scale the Brownian bridge by defining $\tilde{X}^R_t = X_t/R$. The starting and ending positions 
$\tilde{x} = \tilde{X}^R_0, \tilde{y}= \tilde{X}^R_1$ now satisfy $|\tilde{x}|= |\tilde{y}|=1$ and the 
process $\tilde{X}^R_t$ is stopped upon  hitting the ball of radius $1/ 2$. However the process 
$\tilde{X}^R_t$ has reduced variance. Indeed, in law we have the equality
\[
\tilde{X}^R_t = (1-t) \tilde{x} + t \tilde{y} + (B_t - t B_1)/R.
\]
As $R \to \infty$ the process converges to the straight line $\tilde{X}_t = (1-t) \tilde{x} + t \tilde{y}$. For
this limiting process the basic construction is deterministic.
If $\tilde{x} \cdot \tilde{y} \geq 0$ then the straight line never gets closer to the origin
than $2^{-1/2}$. For large $R$ a large deviations estimate shows that deviations away from 
the straight line are exponentially unlikely and (\ref{2.8}) follows. 
To obtain (\ref{2.7}) one again considers the straight line $\tilde{X}_t$  and maximizes 
$\tilde{X}_{\sigma}/(\tau- \sigma)^{1/2}$ over those starting and ending points $\tilde{x}, \tilde{y}$ for
which $\tilde{X}_{\sigma} \cdot \tilde{X}_{\tau} <0$. The maximum occurs, for example, when
$\tilde{X}_{\sigma}/(\tau- \sigma)^{1/2}= (1/2) e_1$ and $ \tilde{X}_{\tau}/(\tau- \sigma)^{1/2}= (1/2) e_2$.
A little trigonometry show that either  $\tilde{X}_{\sigma} \cdot \tilde{X}_{\tau} <0$ or else
$\tilde{X}_{\sigma}/(\tau- \sigma)^{1/2} \leq \tilde{\gamma}$ for some $\tilde{\gamma} \in (0,1)$. 
By taking $\gamma \in (\tilde{\gamma},1)$ a large deviations argument yields (\ref{2.7}).
\qed

\noindent
Applying  (\ref{2.8}) to the bound (\ref{2.6}) we find, when 
$|x|=|y|=R$ and $x \cdot y \geq 0$,
\begin{equation} \label{2.9}
F_K(x,y,1) \leq \exp (4 \eta R^{-2}) \left( 1 + c_5 e^{-c_6 R} \sup_{|x|=|y| \geq R/2} F_K(x,y,1) \right).
\end{equation}

Now we wish to iterate the basic construction to define a Markov chain
$(x(n),y(n))_{n=0,1,\ldots}$ on $\left( \mathbf{R}^{d} \cup \{\Delta\} \right)^2$.
Throughout $|x(n)|=|y(n)|$ or $x(n)=y(n)=\Delta$ will hold.
$\Delta$ is cemetery state from which there is no return.
It will be convenient to set $F(\Delta, \Delta, 1)=F_K(\Delta,\Delta,1) =1$.
We set $x(0)=x, y(0)=y$. Suppose $x(n),y(n)$ have been defined and are not
equal to $\Delta$. Then we repeat the basic construction described above, but started at
the radius $R = |x(n)| = |y(n)| $. We define
\[
\left\{
\begin{array}{ll}
x(n+1) = X_{\sigma}/ (\tau-\sigma)^{1/2}, \; y(n+1) =  X_{\tau}/(\tau-\sigma)^{1/2}
& \mbox{on $ \{ \sigma < \tau\}$,} \\
x(n+1) = y(n+1) =  \Delta
& \mbox{on $\{ |X_t| > R/2, \, \forall t \in [0,1]\}$.} 
\end{array} \right.
\]

We will choose a constant $R_0 \in [1,R]$ shortly. 
Define stopping times for $(x(n),y(n))$ as follows
\begin{eqnarray*}
N_1 & = & \inf \{n: x(n)=y(n)= \Delta\}, \\
N_2 & = & \inf \{n: |x(n)| \leq R_0 \}, \\
N_3 & = & \inf\{n: x(n) \cdot y(n) \geq 0\}, \\
N_4 & = & \inf \{n: |x(n)| > \gamma |x(n-1)| \}.
\end{eqnarray*}
Let $N= N_1 \wedge N_2 \wedge N_3 \wedge N_4$.
Technically we should define $|\Delta|$ to make these times well defined, but we adopt the
convention that if $N \geq k$ and $N_1=k$ then $N_2 = N_3=N_4= \infty$.
Note that $N$ is a bounded stopping time since
if $N_4$ has not occurred then  $N_2 \leq N_0$, where $R \alpha^{N_0} \leq R_0$.
We now expand $F_K(x,y,1)$ as in (\ref{2.6}) to find
\[
F_K(x,y,1) \leq  \sum_{n=1}^{N_0} E \left[ \mathbf{1} (N=n) 
\exp \left( 4 \eta  ( |x(0)|^{-2} + \ldots + |x(n-1)|^{-2}) \right)  
F_K(x(n), y(n),1) \right]. 
\]
On $\{N =n \}$ we know that $R_0 \leq |x(n-1)| \leq \gamma |x(n-2)| \leq \gamma^2 |x(n-3)| \leq \ldots$ 
Hence on this set
\[
\exp \left( 4 \eta  ( |x(0)|^{-2} + \ldots + |x(n-1)|^{-2}) \right)
 \leq \exp \left( \frac{4 \eta}{R_0^2 (1-\gamma^2)} \right).
\]
We choose $R_0 $ large enough that this exponential is bounded by $2$.
This leads to the simpler bound 
\begin{equation} \label{2.10}
F_K(x,y,1) \leq 2 E\left[F_K(x(N), y(N),1) \right] 
\end{equation}
We now find various estimates for  $ E\left[F_K(x(N), y(N),1) \right]$ depending on the value
of $N$.  When $N=N_1$ we have, by definition,
\begin{equation} \label{2.11}
\mathbf{1}(N=N_1)  F_K(x(N), y(N),1) =1.
\end{equation}
When $N=N_2$ we have $|x(n)| \in [R_0/2,R_0]$ and so we can bound
\begin{equation} \label{2.12}
\mathbf{1}(N=N_2)  F_K(x(N), y(N),1) \leq \sup_{|x| = |y| \in [R_0/2, R_0] } F(x,y,1).
\end{equation}
When $N=N_3>N_2$ we have $|x(N)|=|y(N)| \geq R_0$ and $x(N) \cdot y(N) \geq 0$ and we may use 
(\ref{2.9}) to bound $F_K(x(N),y(N),1)$. By choosing $R_0 $ large enough this gives 
the bound 
\begin{equation} \label{2.13}
\mathbf{1}(N=N_3>N_2)  F_K(x(N), y(N),1) \leq 2 \left( 1 + \frac{1}{16} 
 \sup_{|x|=|y| \geq R_0/2} F_K(x,y,1) \right).
\end{equation}
Finally when $N=N_4 < N_2 \wedge N_3 $
we simply bound  
\[
E \left[ \mathbf{1} (N=N_4 < N_2 \wedge N_3) F_K(x(N),y(N),1) \right]
\leq  P(N=N_4 < N_2 \vee N_3) \sup_{|x|=|y| \geq R_0/2}F_K(x,y,1).
\]
We now claim that
\begin{equation} \label{2.14}
\lim_{R \to \infty} \sup_{|x|=|y| \geq R} P(N=N_4< N_2 \wedge N_3) =0.
\end{equation}
Indeed we may apply Lemma \ref{newlemma}
to see that
\[
P(N=N_4 =k+1 <  N_2 \wedge N_3| (x(j),y(j)) \; j=0,1,\ldots,k) \leq c_3 \exp( - c_4 |x(k)|) \mathbf{1}(N>k)
\]
So
\begin{eqnarray*}
\lefteqn{ \sum_{k=1}^{\infty} P(N=N_4=k < N_2 \wedge N_3)}\\
& \leq & \sum_{k=1}^{\infty} E\left[c_3 \exp( - c_4 |x(k-1)|) \mathbf{1}(N>k-1)\right] \\
& \leq & E \left[ c_3 \sum_{k=0}^{N-1} \exp( - c_4 |x(k)|) \right] \\
& \leq & c_3 \sum_{k=0}^{\infty} \exp( - c_4 R_0 \gamma^{-k}) \\
& \to & 0 \quad \mbox{as $R_0 \to \infty$.}
\end{eqnarray*}
Using the claim (\ref{2.14}) we may choose $R_0$ large enough that 
\begin{equation} \label{2.15}
E \left[ \mathbf{1}(N=N_4>N_3 \vee N_2)  F_K(x(N), y(N),1) \right] 
\leq \frac18  \sup_{|x|=|y| \geq R_0/2} F_K(x,y,1).
\end{equation}
Choosing $R_0$ large enough that all four estimates 
(\ref{2.11}), (\ref{2.12}), (\ref{2.13}), (\ref{2.15}) hold, we substitute them into
 (\ref{2.10}) to obtain 
\begin{eqnarray*}
F_K(x,y,1) & \leq & 6 + 2 \sup_{|x|=|y| \in [R_0/2,R_0]} F(x,y,1) + \frac12 \sup_{|x|=|y| \geq R_0/2}
F_K(x,y,1) \\
& \leq & 6 + 3 \sup_{|x|=|y| \in [R_0/2,R_0]} F(x,y,1) + \frac12 \sup_{|x|=|y| \geq R_0}
F_K(x,y,1).
\end{eqnarray*}
Take  the supremum  over $x,y$ in $\{|x|=|y| \geq R_0\}$ of the left hand side
to obtain  
\[
\sup_{ |x|, |y| \geq R_0/2} F_K(x,y,1) \leq  
12 + 6 \sup_{|x|=|y| \in [R_0/2,R_0]} F(x,y,1).
\]
Letting $K \to \infty$ we have bounded $F(x,y,1)$ on the set
$\{ |x|=|y| \geq R_0/2 \}$. Together with (\ref{2.5})
this completes the proof of the main estimate.

\noindent 
\textbf{Remarks}

\textbf{1.} The moment $ E_{a,b,t}^{(d)} \left[ \exp ( \eta \int_{0}^{t} R_s^{-2} ds)  \right]$
is infinite for $\eta > (d-2)^2/8$. This follows since the formula (\ref{2.2})
cannot be  analytically extended, as a function of  $\eta$,  into the  region $\{z: Re(z) < r\}$ 
for any $r > (d-2)^2/8$. This strongly suggests there are no solutions to (\ref{1.1}) having 
finite second moments $E[(u_t(f))^2]$ when $ \kappa > (d-2)/2$. Similarly, the blow-up of the 
Brownian exponential moment suggests
there should be no solutions to (\ref{1.1}) with finite second moments 
for any $ \kappa>0$ when the noise
has covariance (\ref{1.2})  with $p>2$.  

\textbf{2.} As indicated in Subsection \ref{Tools}, higher moments are controlled by the 
Brownian exponential moments (\ref{1.13}). Using H\"{o}lder's inequality we find 
\begin{eqnarray*}
\lefteqn{E_{0,x_1,\dots,x_n}^{t,y_1,\dots,y_n}
\left[\exp\left(\sum_{1\leq j<k \leq n}\int_{0}^{t} \frac{\kappa^2} {\left|X^{(j)}_s-X^{(k)}_s \right|^{2}}
 ds \right) \right]}\\
& \leq & \prod_{1 \leq j < k \leq n} \left( E_{0,x_j,x_k}^{t,y_j,y_k}
\left[\exp\left( \int_{0}^{t} \frac{n(n-1) \kappa^2} {\left|X^{(j)}_s-X^{(k)}_s \right|^{2}}
 ds \right) \right] \right)^{1/n(n-1)}.
\end{eqnarray*}
The exponential moment calculated in this section shows that this  is finite 
when $n(n-1) \kappa^2/2 \leq (d-2)^2/8$. This should lead to the solutions
to (\ref{1.1}) having finite moments $E[(u_t(f))^n]$ when $\kappa \leq (d-2) (4n(n-1))^{-1/2}$.
We do not think this simple H\"{o}lder argument leads to the correct critical values
for the existence of higher moments.
%
%
\section{Existence of Solutions}
%
%
\label{s4} \setcounter{equation}{0} 

In this section we give a construction of solutions to (\ref{1.1}) using the chaos expansion 
(\ref{1.15}). However, it is hard to show from the series expansion that the resulting 
solution is a non-negative measure. For that purpose we give a second construction as a limit of 
less singular SPDEs. A comparison theorem will show that the approximating equations have 
solutions which are non-negative functions implying that the limit must also be non-negative.
Finally we show that the two constructions yield the same process and that it is a  solution 
of (\ref{1.1}).

We first construct a noise $F$ with the desired covariance. Let 
$g(x)=c_7 |x|^{-(d+2)/2}$. A simple calculation shows, for a suitable value of the constant $c_7$, 
that the convolution $g * g(z) = |z|^{-2}$. Now let $W$ be an adapted space-time white noise on 
$\mathbf{R}^d \times [0,\infty)$ on some filtered probability space 
$(\Omega, \mathcal{F}, \{ \mathcal{F}_t\} ,P)$. Define, for $f: \mathbf{R}^d \rightarrow \mathbf{R}$ 
that is bounded, measurable and of compact support, 
\begin{equation}  \label{3.1}
F(t,f) = \int^t_0 \int (f * g)(z) W(dz,ds).
\end{equation}
It is straightforward to show that $F(t,f)$ is well defined, is a Gaussian martingale, and that 
\[
\langle F(\cdot,f)\rangle_t = t \; \int \! \int  \frac{f(y) f(z)}{|y-z|^2} dy \, dz. 
\]
If we write $F(t,A)$ when $f=I_A $ then $\{F(t,A):t \geq 0, A \subseteq  \mathbf{R}^d\}$ is a 
martingale measure and hence (see \cite{Walsh86} Chapter 2) can be used to define a stochastic 
integral  $\int^t_0 \int h(s,y) F(dy,ds) $ for suitable predictable integrands $h$ so that 
\[
\left[ \int_0^{\cdot} \! \int h(s,y)F(dy,ds) \right]_t = \int^t_0 \! 
\int \!\int \frac{h(s,y) h(s,z)}{|y-z|^2} dy \, dz \, ds. 
\]

Next we show that the expansion (\ref{1.15}) converges.
\begin{Lemma} \label{series-converges} 
Suppose $\mu \in \mathcal{H}_{\alpha}$. Then, for 
$f \in \mathcal{C}_c$, the series 
$\sum_{n=0}^{\infty}I^{(n)}_t(f,\mu )$, defined by (\ref{1.16}) 
and (\ref{1.17}), converges in $L^2$. Moreover 
\begin{eqnarray}  
E\left[ \left( \sum_{n=0}^{\infty} I^{(n)}_t(f, \mu ) \right)^2 \right] 
& = & \sum_{n=0}^{\infty} E \left[ \left( I^{(n)}_t(f, \mu ) \right)^2 \right]  \nonumber \\
&=& \int_{\mathbf{R}^{4d}} f(y')f(x') G_t(x-x')G_t(y-y') \label{3.2} \\
&& \hspace{.4in} \cdot E_{0,x,y}^{t,x',y'}\left[\exp\left( \int_{0}^{t} \frac{\kappa^2}{\left|X_s^1-X_s^2\right|^2}
ds\right)\right] \; dx' dy' \mu(dx) \mu(dy). \nonumber
\end{eqnarray}
\end{Lemma}

\noindent
\textbf{Proof.} We first check that the right hand side of (\ref{3.2}) is finite. Using 
the fact that $(X^1_t-X^t_2)/\sqrt{2}$ is a Brownian bridge from $x-y$ to $x'-y'$ we may use 
Lemma   \ref{exponentialmoment} to obtain 
\begin{eqnarray*}
\lefteqn{ \int_{\mathbf{R}^{4d}} f(y')f(x') G_t(x-x')G_t(y-y') E_{0,x,y}^{t,x',y'}
\left[\exp\left( \int_{0}^{t} \frac{\kappa^2}{\left|X_s^1-X_s^2\right|^2}
ds\right)\right] \mu(dx) \mu(dy) dx' dy' } \\
& \leq & C(t) \int_{\mathbf{R}^{4d}} f(y')f(x') G_t(x-x')G_t(y-y') 
\left(1 + |x-y|^{-\alpha} |x'-y'|^{-\alpha} \right) \mu(dx) \mu(dy) dx' dy' .  
\end{eqnarray*}
Now estimates as in Lemma \ref{momentcor1} show this expression is finite.

The multiple Wiener integrals of different orders are orthogonal, if they have finite second moments; 
that is if $m\ne n$, and if 
\begin{equation}  \label{3.3}
E\left[ \left( I^{(k)}_t(f, \mu) \right)^2 \right] < \infty
\end{equation}
for $k=m,n$, then 
\[
E\left[I^{(m)}_t(f, \mu)I^{(n)}_t(f, \mu)\right]=0. 
\]
It is therefore enough to establish the second equality in (\ref{3.2}) since this implies 
(\ref{3.3}), and then orthogonality of the terms in the series implies the first equality in 
(\ref{3.2}). First note that, with $s_{n+1}=t$, 
\begin{eqnarray} 
E\left[ \left( I^{(n)}_t(f) \right)^2 \right] &=& \kappa^{2n} \int_{0}^{t} 
\int_{0}^{s_n}{\ldots}\int_{0}^{s_2}
 \int_{\mathbf{R}^{2(n+1)d}} f(y_{n+1})f(z_{n+1})dy_{n+1}dz_{n+1} 
G_{s_1} \mu (y_1) G_{s_1} \mu (z_1) \nonumber \\
&& \hspace{.2in} \cdot \prod_{i=1}^{n}\left[G_{s_{i+1}-s_{i}}(y_{i+1}-y_{i})
G_{s_{i+1}-s_{i}}(z_{i+1}-z_{i})|y_i-z_i|^{-2}dy_idz_ids_i\right]. \label{3.4}
\end{eqnarray}
Expanding  the exponential in the final term of (\ref{3.2}) we have 
\begin{eqnarray*}
E_{0,x,y}^{t,x',y'}\left[\exp\left( \int_{0}^{t} \frac{\kappa^2}{\left|X_s^1-X_s^2\right|^2} 
ds \right) \right] 
&=& 1+ E_{0,x,y}^{t,x',y'}\left[\sum_{n=1}^{\infty}\frac{1}{n!} \int_{0}^{t}{\ldots}
\int_{0}^{t} \prod_{i=1}^{n} \frac{\kappa^2 ds_i}{\left| X_{s_i}^1-X_{s_i}^2 \right|^2} \right] \\
&=& 1+  \sum_{n=1}^{\infty}\int_{0}^{t}\int_{0}^{s_n}{\ldots}\int_{0}^{s_2}
E_{0,x,y}^{t,x',y'}\left[\prod_{i=1}^{n} 
\frac{\kappa^2 ds_i}{\left|X_{s_i}^1-X_{s_i}^2\right|^2}\right].
\end{eqnarray*}
Substituting this sum into the right hand side of  (\ref{3.2}), one may match,  by
using the finite dimensional distributions of the Brownian bridge, the nth term with the 
expression $E[I^{(n)}_t(f)]^2$ in (\ref{3.4}).  \qed

The chaos expansion defines a linear random functional on test functions (in that there is a 
possible null set for each linear relation). Also this linear random functional satisfies the
moment bounds (\ref{1.6}) and (\ref{1.7}).
The second moment bound implies that there is a regularization 
(see \cite{Ito84} Theorem 2.3.3),  ensuring there is a random distribution $u_t$ so that
\begin{equation}  \label{3.5}
u_t(f) = \sum_{i=0}^{\infty} I^{(n)}_t(f, \mu) \quad \mbox{for all $f
\in \mathcal{C}_c$, almost surely.}
\end{equation}
To show that $u_t$ is actually a random measure we 
now construct a sequence of SPDE approximations to (\ref{1.1}). We will index our 
approximations by numbers $\varepsilon>0$. Recall that $h(x) :=|x|^{-2}=(g*g)(x)$, where 
$g(x)=c_7 |x|^{-(d+2)/2}$. Let 
\[
g^{(\varepsilon)}(x)=\left(c_7 |x|^{-(d+2)/2}\right) \wedge \varepsilon^{-1}, \quad 
\mbox{and} \quad h^{(\varepsilon)}(x) = \left( g^{(\varepsilon)}*g^{(\varepsilon)}\right) (x). 
\]
As $\varepsilon \downarrow 0$ we have $g^{(\varepsilon)}(x) \uparrow g(x)$ and 
$h^{(\varepsilon)}(x)\uparrow h(x).$ We can construct, as in (\ref{3.1}), a mean zero 
Gaussian field $F^{(\varepsilon)}(t,x)$ with covariance 
\[ 
E\left[\dot F^{(\varepsilon)}(t,x)\dot F^{(\varepsilon)}(s,y)\right] =
\delta(t-s)h^{(\varepsilon)}(x-y). 
\]
We consider the approximating SPDE 
\begin{equation}  \label{3.6}
\frac{\partial u^{(\varepsilon)} }{\partial t} =  \frac12 \Delta u^{(\varepsilon)} +
\kappa u^{(\varepsilon)} \dot{F^{(\varepsilon)} }, \quad u_0^{(\varepsilon)} = \mu^{(\delta)},
\end{equation}
with the initial condition $\mu^{(\delta)}= G_{\delta} \mu $, for some
$\delta = \delta(\varepsilon)>0$ to be chosen later. Since the correlation is 
continuous in $x$ and $y$, standard results give existence and uniqueness of a 
non-negative, continuous, function-valued solution $u^{(\varepsilon)}_t(t,x)$. 
Moreover we may represent the solutions in terms of a chaos expansion 
\[
u^{(\varepsilon)}_t(f)=\sum_{n=0}^{\infty} I^{(n,\varepsilon)}_t(f, \mu^{(\delta)}) 
\]
where the terms $I^{(n,\varepsilon)}_t(f, \mu^{(\delta)}) $ are defined as in (\ref{1.16}) 
and (\ref{1.17}) except that $\mu,F$ are replaced by $\mu^{(\delta)},F^{(\varepsilon)}$.
We now connect the approximations with the original series construction.
\begin{Lemma} \label{ell-two} 
Suppose that $\mu \in \mathcal{H}_{\alpha+}$. Then we may define $I^{(n)}_t(f,\mu)$ and  
$I^{(n,\varepsilon)}_t(f, \mu^{(\delta)})$ on the same probability space so that, for 
suitably chosen $\delta(\varepsilon) >0$, fixed $t \geq 0$ and $f \in \mathcal{C}_c$, 
\[
u^{(\varepsilon)}_t(f) \rightarrow  
u_t(f)  \quad \mbox{in $L^2$ as $\varepsilon \rightarrow 0$.} 
\]
Hence the chaos expansion  (\ref{3.5}) defines a random measure $u_t(dx)$, for each 
$\mu \in \mathcal{H}_{\alpha+}$ and $t \geq 0$. 
\end{Lemma}
\noindent
\textbf{Proof.} Let $\dot W(t,x)$ be a space-time white noise on $[0,\infty)\times\mathbf{R}^d$ 
and construct both the noises $F$ and $F^{(\varepsilon)}$ using $W$ as in (\ref{3.1}). Using 
the convergence of both the $F$ and $F^{(\varepsilon)}$ chaos expansions and the 
orthogonality  of multiple Wiener integrals of different orders, we find 
\begin{eqnarray} 
\lefteqn{E\left[ \left( u^{(\varepsilon)}_t(f) - u_t(f) \right)^2 \right]}
\nonumber  \\
& = & \sum_{n=0}^{\infty} E\left[ \left( I^{(n)}_t(f, \mu )- 
I^{(n,\varepsilon)}_t(f, \mu^{(\delta)}) \right)^2 \right]  \label{3.7} \\
& \leq & 2 \sum_{n=0}^{\infty} E\left[ \left( I^{(n)}_t(f, \mu )- I^{(n, \varepsilon)}_t(f, \mu) 
\right)^2 \right] + 2 \sum_{n=0}^{\infty} E\left[ \left( I^{(n, \varepsilon)}_t(f, \mu ) -
I^{(n,\varepsilon)}_t(f, \mu^{(\delta)}) \right)^2 \right].  \nonumber
\end{eqnarray}
We show separately that both sums on the right hand side of (\ref{3.7}) converge to zero
as $\varepsilon \downarrow 0$. We use the telescoping expansion, for $n \geq 1$, 
\begin{eqnarray*}
\lefteqn{I^{(n)}_t(f, \mu)- I^{(n,\varepsilon)}_t(f, \mu)} \\
&=& \kappa^n \int_{0}^{t}{\ldots}\int_{0}^{s_2} \int_{\mathbf{R}^{(n+1)d}} G_{s_1} \mu 
(y_1)f(y_{n+1})dy_{n+1} \prod_{i=1}^{n}\left[G_{s_{i+1}-s_{i}}(y_{i+1}-y_{i}) F(dy_i, ds_i)\right] \\
&& -\kappa^n \int_{0}^{t}{\ldots}\int_{0}^{s_2} \int_{\mathbf{R}^{(n+1)d}}
G_{s_1} \mu(y_1)f(y_{n+1})dy_{n+1} \prod_{i=1}^{n} \left[G_{s_{i+1}-s_{i}}(y_{i+1}-y_{i}) 
F^{(\varepsilon)}(dy_i, ds_i)\right] \\
& = & \sum_{m=1}^{n} J^{(n,m,\varepsilon)}_t(f,\mu)
\end{eqnarray*}
where $J^{(n,m,\varepsilon)}_t(f) $ is defined to equal 
\begin{eqnarray*}
&& \kappa^n \int_{0}^{t}{\ldots} \int_{0}^{s_2} \int_{\mathbf{R}^{(n+1)d}}
G_{s_1}\mu(y_1)f(y_{n+1})dy_{n+1} \prod_{i=1}^{n}G_{s_{i+1}-s_{i}}(y_{i+1}-y_{i}) \\
&& \hspace{.3in} \cdot \prod_{i=1}^{m-1}F(dy_i, ds_i) \cdot  \left[F(dy_m, ds_m) - 
F^{(\varepsilon)}(dy_m, ds_m)\right] \cdot \prod_{i=m+1}^{n}F^{(\varepsilon)}(dy_i, ds_i)
\end{eqnarray*}
and where a product over the empty set is defined to be 1. The isometry for the stochastic integral gives 
\begin{eqnarray*}
\lefteqn{E\left[ \left(J^{(n,m,\varepsilon)}_t(f, \mu)\right)^2 \right]} \\
&=& \kappa^{2n}\int_{0}^{t}{\ldots}\int_{0}^{s_2} \int_{\mathbf{R}^{2(n+1)d}} G_{s_1} 
\mu(y_1) G_{s_1} \mu(z_0) f(y_{n+1})f(z_{n+1}) dy_{n+1}dz_{n+1} \\
&&  \hspace{.3in}  \cdot  \prod_{i=1}^{n} \left[G_{s_{i+1}-s_{i}}(y_{i+1}-y_{i})
G_{s_{i+1}-s_{i}}(z_{i+1}-z_{i})dy_idz_ids_i\right] \\
&&  \hspace{.6in} \cdot  \prod_{i=1}^{m-1}h(y_i-z_i) 
\left[ \left(g-g^{(\varepsilon)}\right) *\left(g-g^{(\varepsilon)}\right) \right]  (y_m-z_m)
\prod_{i=m+1}^{n}h^{(\varepsilon)}(y_i-z_i).
\end{eqnarray*}
Note that $0 \leq [(g-g^{(\varepsilon)}) *(g-g^{(\varepsilon)})](x) \leq h(x)$ and that 
$[(g-g^{(\varepsilon)}) *(g-g^{(\varepsilon)})](x) \downarrow 0$ as $\varepsilon \to 0$.
Using  the finiteness of  $E \left[(I_t^{(n)}(f,\mu))^2 \right]$, the dominated convergence theorem 
implies that  $E\left[ \left(J^{(n,m,\varepsilon)}_t(f,\mu) \right)^2 \right] \downarrow 0$ and therefore 
\begin{equation}  \label{3.8}
\lim_{\varepsilon\downarrow 0} E\left[ \left( I^{(n,\varepsilon)}_t(f,\mu) - 
I^{(n)}_t(f,\mu) \right)^2 \right] = 0.
\end{equation}
The isometry, and $h^{(\varepsilon)}(x) \leq h(x)$, imply that 
\begin{eqnarray}
E \left[ \left( I^{(n)}_t(f) -I^{(n,\varepsilon)}_t(f)\right)^2 \right]
& \leq & 2 E\left[ \left( I^{(n)}_t(f) \right)^2 \right]
+ 2 E\left[ \left(I^{(n,\varepsilon)}_t(f)\right)^2\right]  \nonumber \\
\label{3.9} & \leq & 4 E\left[ \left(I^{(n)}_t(f)\right)^2 \right]. 
\end{eqnarray}
Now, using (\ref{3.8}), (\ref{3.9}), the convergence of the series
$\sum_{n=0}^{\infty}E\left[\left( I^{(n)}_t(f)\right)^2 \right]$, and the dominated convergence theorem,
the first term on the right hand  side of (\ref{3.7}) goes to zero as
$\varepsilon \downarrow 0$.

We now show that for fixed $\varepsilon>0$ the second term on the right hand side  of
(\ref{3.7}) converges to zero as $ \delta \downarrow 0$. Recall the initial
condition was $\mu^{(\delta)}= G_{\delta} \mu$ for some $\delta= \delta(\varepsilon)>0$. 
But for fixed $\varepsilon$ the isometry shows, as in Lemma \ref{series-converges}, that
\begin{eqnarray*} 
\lefteqn{ \sum_{n=0}^{\infty} E\left[ \left( I^{(n, \varepsilon)}_t(f, \mu ) 
- I^{(n, \varepsilon)}_t(f, \mu^{(\varepsilon)} ) \right)^2 \right]} \\
&=& \int_{\mathbf{R}^{4d}} f(y')f(x') G_t(x-x')G_t(y-y') \\
&& \hspace{.4in} \cdot E_{0,x,y}^{t,x',y'}\left[\exp\left( \int_{0}^{t} \kappa^2 h^{(\varepsilon)}(X_s^1-X_s^2) ds\right)\right]  \; dx' dy' \left( \mu - \mu^{(\delta)} \right) (dx) \left( \mu - \mu^{(\delta)} \right) (dy).  \end{eqnarray*}
When $\varepsilon>0$,  the Brownian bridge expectation is a bounded continuous function
of $x,y,x',y'$ and the convergence to zero as $\delta \downarrow 0$ is clear. This completes the 
proof of the $L^2$ convergence stated in the lemma.

The $L^2$ boundedness of $u_t^{(\varepsilon)}(f)$, for each $f \in \mathcal{C}_c$, implies that 
$\{u^{(\varepsilon)}_t(x) dx\}$ is a tight family of random Radon measures. 
The $L^2$ convergence of $u_t^{(\varepsilon)}(f)$ implies that there is a random measure
$u_t$ satisfying (\ref{3.5}) and that $u_t^{(\varepsilon)} \to u_t$ in distribution
as $\varepsilon \to 0$. \qed

It remains to show that $\{u_t(dx)\}$ is a solution of (\ref{1.1}), and for this we must show that
there is a continuous version of the process $t \rightarrow u_t$ and that it satisfies the 
martingale problem (\ref{1.4}) and (\ref{1.5}).
Fix $f \in \mathcal{C}_c^2$. From the definition (\ref{1.17}) we have, for $n \geq 1$,
\[ 
I^{(n)}_{s}(y,z) = \kappa \int^{s}_0 \int G_{s-r}(y-y') I^{(n-1)}_{r} (y',z) F(dy',dr).
\]
Then using a stochastic Fubini theorem (see \cite{Walsh86} Theorem 2.6), and the fact 
that $G_t * f(y)$ solves the heat equation,  we have, for $n \geq 1$,
\begin{eqnarray*}
\int^t_0  I^{(n)}_s ( \frac12 \Delta f,\mu) ds & = & \frac12  \int^t_0  \int \! \int I^{(n)}_s (y,z)
\Delta f(y) \mu(dz) dy ds \\
& = & \frac{\kappa}{2} \int^t_0 \int \! \int \left(  \int^s_0  \int G_{s-r}(y-y')  
I^{(n-1)}_{r} (y',z) F(dy',dr) \right) \Delta f(y) \mu(dz) dy ds  \\
& = & \frac{\kappa}{2} \int^t_0 \int \left( \int^t_r  G_{s-r} * \Delta f(y') ds \right) \int I^{(n-1)}_{r} (y',z) 
\mu(dz)   F(dy',dr)  \\
& = & \kappa \int^t_0 \int \left( G_{t-r} * f(y') - f(y') \right) \int I^{(n-1)}_{r} (y',z) \mu(dz) F(dy',dr)  \\
& = & I^{(n)}_t(f,\mu) - \kappa \int^t_0 \int \! \int  f(y')  I^{(n-1)}_{r} (y',z) \mu(dz) F(dy',dr). 
\end{eqnarray*}
Rearranging the terms, we see that for each $n \geq 1$, the process
\begin{equation} \label{3.10}
z^{(n)}_t(f) := I^{(n)}_t(f, \mu) - \int^t_0 I^{(n)}_s (\frac12 \Delta f,\mu) ds 
=  \kappa \int^t_0 \int \! \int f(y)  I^{(n-1)}_{s} (y,z) \mu(dz) F(dy,ds) 
\end{equation}
is a continuous martingale.  We now define
\[
u_{N,t}(x) = \sum_{k=0}^N I^{(k)}(x,x') \mu(dx'), \quad z_{N,t}(f) = \sum_{k=1}^N z^{(k)}_t(f).
\]
Then, for $f \in \mathcal{C}^2_c$ and $N \geq 1$,
\begin{equation} \label{3.11}
u_{N,t}(f) = \mu(f) + \int^t_0 u_{N,s}(\frac12 \Delta f)ds + z_{N,t}(f).
\end{equation}
Lemma  \ref{series-converges} implies that $E[ (u_{N,t}(\Delta f) - u_t(\Delta f))^2 ]$ converges
monotonically to zero. Using the domination from Lemma \ref{momentcor1} part i) we have
\[
E\left[ \sup_{t \leq T} \left| \int^t_0 u_{N,s}(\frac12 \Delta f)ds 
- \int^t_0 u_s(\frac12 \Delta f) ds \right|^2 \right] \to 0. 
\]
Lemma  \ref{series-converges} also implies that $ z_{N,t}(f)$ converges in $L^2$ to $z_t(f)$
and by Doob's inequality
\[
E\left[ \sup_{t \leq T} \left| z_{N,t}(f) - z_t(f) \right|^2 \right] \to 0.
\]
This uniform convergence and (\ref{3.11}) shows  there is a continuous version of both
$t \to z_t(f)$ and $t \to u_t(f)$.  Using this for a suitable countable class of $\mathcal{C}_c^2$ test functions $f$ shows that there is a continuous version (in the vague topology) of $t \rightarrow u_t$.

Now we calculate the quadratic variation $z_t(f)$, which is the 
$\mathbf{L}^1$ limit of $\left\langle z_{N,\cdot}(f) \right\rangle_t$.
It is enough to consider the case  $f \geq 0$.  Using (\ref{3.10}) we have
\begin{eqnarray}
\left\langle z_{N+1,\cdot}(f) \right\rangle_t  & = & \sum_{k=1}^{N+1} \sum_{l=1}^{N+1} 
\kappa^2 \int^t_0 \int_{\mathbf{R}^{4d}} \frac{f(x) f(y)}{|x-y|^2} I^{(k-1)}_s (x,x') 
I^{(l-1)}_s (y,y') \mu(dx') \mu(dy') dx \, dy \, ds \nonumber \\
& = & \kappa^2 \int^t_0 \int_{\mathbf{R}^{2d}} \frac{f(x) f(y)}{|x-y|^2}u_{N,s}(x)
u_{N,s}(y) dx \, dy \, ds \nonumber \\
& \to & \kappa^2 \int^t_0 \int_{\mathbf{R}^{2d}} \frac{f(x)f(y)}{|x-y|^2}
u_s(dx) u_s(dy) ds. \label{3.12}
\end{eqnarray}
We need to justify this final convergence and we split the task into two terms
\begin{eqnarray*}
\lefteqn{\int^t_0 \int_{\mathbf{R}^{2d}} \frac{f(x) f(y)}{|x-y|^2}u_{N,s}(x)dx \, 
(u_{N,s}(y)dy - u_s(dy)) \, ds} \\
& & + 
 \int^t_0 \int_{\mathbf{R}^{2d}} \frac{f(x)f(y)}{|x-y|^2}
(u_{N,s}(x) dx - u_s(dx))  u_s(dy) ds.
\end{eqnarray*}
We show the first term converges to zero in $L^1$, the argument for the second 
term is the same. We use the fact that $|x-y|^{-2}$ is a convolution of $c_7 |z|^{-(d+2)/2}$ with 
itself to see that
\begin{eqnarray*}
\lefteqn{\int_{\mathbf{R}^{2d}} \frac{f(x) f(y)}{|x-y|^2}u_{N,s}(x)dx \,  (u_{N,s}(y)dy - u_s(dy))} \\
& = & c_7 \int_{\mathbf{R}^{3d}} \frac{f(x) f(y)}{|x-z|^{(d+2)/2} |y-z|^{(d+2)/2}}u_{N,s}(x)dx
(u_{N,s}(y)dy - u_s(dy)) dz \\
& = & c_7 \int_{\mathbf{R}^d}  u_{N,s}(f_z) (u_{N,s}(f_z) - u_s(f_z)) dz
\end{eqnarray*}
 where $f_z(x) = f(x) |x-z|^{-(d+2)/2}$. Hence, by the Cauchy-Schwartz inequality,
\begin{eqnarray}
\lefteqn{ E \left[ \left| \int^t_0 \int_{\mathbf{R}^{2d}} \frac{f(x) f(y)}{|x-y|^2}u_{N,s}(x)dx
(u_{N,s}(y)dy - u_s(dy)) \, ds \right| \right]} \label{3.13} \\ 
& \leq & c_7 \left( \int^t_0 \int  E\left[ (u_{N,s}(f_z))^2 \right] dz ds \right)^{1/2}
 \left( \int^t_0 \int  E\left[ (u_{N,s}(f_z)-u_s(f_z))^2 \right] dz ds \right)^{1/2} 
\nonumber
\end{eqnarray}
The argument from Lemma \ref{series-converges} shows that $ E\left[ (u_{N,s}(f_z))^2 \right]$
can be bounded uniformly in $N$ by
\[
E\left[ (u_{N,s}(f_z))^2 \right] \leq \int_{\mathbf{R}^{4d}} \! f_z(x') f_z(y') G_s(x-x') G_s(y-y')
\! \left( 1 + \frac{s^{\alpha}}{|x-y|^{\alpha} |x'-y'|^{\alpha}} \right) \! \mu(dx) \mu(dy) dx' dy'.
\]
The same bound holds for $ E\left[ (u_{s}(f_z))^2 \right]$. It is straightforward but lengthy to
estimate this term. We show how to deal with the most singular term only. 
The method is to estimate the $dx' dy' $ integral first using 
the inequalities (\ref{1.9}) and (\ref{1.10}). Applying  H\"{o}lder's 
inequality in the same way as in Lemma \ref{momentcor1}, these inequalities imply that
\begin{eqnarray*}
\lefteqn{ \int_{\mathbf{R}^{2d}} G_s(x-x') G_s(y-y') \frac{f(x') f(y')}{|x'-z|^{(d+2)/2} |y'-z|^{(d+2)/2}
|x'-y'|^{\alpha}} dx' dy'} \\
& \leq & C(a) e^{-a|x|-a|y|} s^{-\alpha/2}
\left( |x-z|^{-(d+2)/2} \wedge s^{-(d+2)/4} \right) \left( |y-z|^{-(d+2)/2} \wedge s^{-(d+2)/4} \right) \\
& \leq & C(\beta,a) e^{-a|x|-a|y|}  s^{-1 -\alpha + (\beta/2)} 
 |x-z|^{-(d+\beta-\alpha)/2} |y-z|^{-(d+\beta-\alpha)/2}. 
\end{eqnarray*}
where we have chosen $\beta \in (\alpha, \alpha+2) $ and $a$ so that 
$\mu \in \mathcal{H}_{\beta}^a$.
To apply Holder's inequality here, splitting the three factors $f(x')f(y')$, 
$ |x'-z|^{-(d+2)/2} |y'-z|^{-(d+2)/2}$ and $|x'-y'|^{\alpha}$, we needed
the bound $\alpha + ((d+2)/2) <d$, 
which is implied by our assumption that $\alpha < (d-2)/2$. 
Substituting this estimate into (\ref{3.13}) we find
\begin{eqnarray*}
\lefteqn{\int^t_0 \int  E\left[ (u_{N,s}(f_z))^2 \right] dz ds} \\
& \leq & C(\beta,a)  \int^t_0 \int_{\mathbf{R}^{3d}} e^{-a|x|-a|y|}  s^{-1 + (\beta/2)} 
|x-y|^{-\alpha} \\
&& \hspace{.3in} \cdot  |x-z|^{-(d+\beta-\alpha)/2} |y-z|^{-(d+\beta-\alpha)/2}  \mu(dx) \mu(dy) dz ds \\
& = & C(\beta,a)  \int^t_0 \int_{\mathbf{R}^{2d}} e^{-a|x|-a|y|}  s^{-1 + (\beta/2)} 
|x-y|^{-\beta} \mu(dx) \mu(dy) dz ds 
\end{eqnarray*} 
which is finite since $\mu \in \mathcal{H}^a_{\beta}$. This bound also gives the domination 
required to see that $\int^t_0 \int  E\left[ (u_{N,s}(f_z)-u_s(f_z))^2 \right] dz ds \to 0$
as $n \to \infty$. 
This finishes the justification of  the convergence in (\ref{3.12}),  identifying  
the quadratic variation  $\langle z_{\cdot}(f) \rangle_t$, and completes the construction 
of a solution $\{u_t(dx)\}$ to (\ref{1.1}) started at $\mu$. 
%
%
\section{Self Duality and Uniqueness}
%
%
\label{s5} \setcounter{equation}{0} 

In this section we establish the self duality of solutions in the following form:
\begin{Proposition} \label{dualityrelation} 
Suppose $\{u_t(dx)\}$ and $\{v_t(dx)\}$ are solutions of (\ref{1.1}), with deterministic initial conditions 
$u_0(dx) =f(x)dx$ and $v_0(dx) = g(x)dx $.  Suppose also that $sup_{x} e^{-a|x|} f(x) < \infty$ for some $a$ 
and that $g(x)$ is bounded and has compact support. Then $u_t(g)$ has the same distribution as  $v_t(f)$.
\end{Proposition}
\noindent
\textbf{Remarks}

\textbf{1}.  The duality formula is immediately clear for the solutions constructed using the chaos expansion 
in Section, \ref{s4} since the expression (\ref{1.17}) for the nth order of the expansion is symmetric 
under the interchange of $y$ and $z$. 
We will show in this section that the self duality relation holds for any solution to (\ref{1.1}). 
We then use the self duality relation to show uniqueness in law for solutions. 

\textbf{2}. Even when working with the martingale problems the self duality relation is heuristically 
clear, as can be seen by applying the technique of Markov process duality  (see Ethier and Kurtz 
\cite{Ethier+Kurtz86} chapter 4). Take $\{u_t(dx)\}$ and $\{v_t(dx)\}$ to be independent solutions to 
(\ref{1.1}). Suppose 
(falsely) that the solutions are function valued and have suitable behavior at infinity such that 
the integrals $u_s (v_{t-s}) $ and $v_{t-s}(u_s)$ are finite and equal by integration by parts. Take a twice 
differentiable $h: [0,\infty) \rightarrow \mathbf{R}$. Applying Ito's formula formally, using the martingale 
problem (\ref{1.4}),  leads to 
\[
\frac{d}{ds} h(u_s ( v_{t-s})) =  (1/2) h^{\prime}(u_s(v_{t-s})) 
\left( u_s(\Delta v_{t-s}) - v_{t-s}(\Delta u_s) \right) + \mbox{martingale terms}. 
\]
Here we have used the cancellation of the two second derivative terms involving $h^{\prime\prime}$ 
after applying Ito's formula for $u_s$ and for $v_{t-s}$. Applying integration by parts the term 
$\left( u_s(\Delta v_{t-s}) - v_{t-s}(\Delta u_s) \right) $ vanishes and this leaves only martingale
terms. Taking expectations and integrating over $s \in [0,t]$ leads to 
\begin{equation}  \label{4.1}
E \left[h( u_t(g) ) \right] =E \left[ h( v_t(f) ) \right]
\end{equation}
which implies the self duality. To make this argument rigorous we shall 
argue using a smoother approximate  duality relation.

\textbf{3}. The self duality relation can be extended to hold for more general initial conditions and to be 
symmetric in the requirements on the initial conditions $\mu$ and $\nu$, as would be expected by the 
symmetry of the chaos expansion. One needs to define certain collision integrals 
$(\mu, \nu)$ between measures in $\mathcal{H}_{\alpha+}$. For example, suppose $\mu, \nu \in 
\mathcal{H}_{\alpha+}$ and for simplicity suppose both are supported in the ball $B(0,R)$. Define 
$f_{\varepsilon}(x) = \int \phi_{\varepsilon}(x-y) \nu(dy)$, the density of the measure 
$\phi_{\varepsilon} * \nu$. Then, if $\{u_t(dx)\}$ is a solutions started at $\mu$, we claim that the 
variables 
\[
u_t(f_{\varepsilon}) = \int \! \int \phi_{\varepsilon}(x-y) u_t (dx) \nu(dy)
\]
are Cauchy in $L^2$ as $\varepsilon \to 0$. Indeed, using the second moment formula 
(\ref{1.7}), a short calculation leads to 
\begin{eqnarray*}
\lefteqn{ E \left[ \left( u_t(f_{\varepsilon}) - u_t(f_{\varepsilon'}) \right)^2 \right]} \\
& = & E \left[ \left( u_t( f_{\varepsilon} - f_{\varepsilon^{\prime}} ) \right)^2 \right] \\
& \leq & C(t,R,\mu) \int \! \int (f_{\varepsilon}(x) - f_{\varepsilon^{\prime}}(x))
 (f_{\varepsilon}(y) - f_{\varepsilon^{\prime}}(y)) (1 + |x-y|^{-\alpha}) dx dy \\
& = & C(t,R,\mu) \| \phi_{\varepsilon} * \nu - \phi_{\varepsilon^{\prime}} * \nu \|_{\alpha}^2.
\end{eqnarray*}
Here we are extending the use of the norm $\| \mu \|_{\alpha}$ to signed measures. Now it is not 
difficult to show that $\| \phi_{\varepsilon} * \nu - \nu \|_{\alpha} \to 0$ as $\varepsilon \to 0$ which 
completes the proof of the Cauchy property. Denoting the $L^2$ limit as $u_t(\nu)$, and constructing 
$v_t(\mu)$ analogously, the duality relation then holds in this extended setting when $\mu,\nu \in 
\mathcal{H}_{\alpha+}^0$, although we make no 
use of it in this paper.
In the rest of this section we give the proof of Proposition \ref{dualityrelation} and deduce uniqueness 
in law and the Markov property. The proof follows from two lemmas, 
the first of which is  an approximate duality relation where we smooth  the measure valued solutions.
\begin{Lemma} \label{doubletestfunction} 
Suppose $\{u_t(dx)\}$ is a solution of (\ref{1.1}) with initial condition $\mu$
and $\{v_t(dx)\}$ is an independent solution with a compactly supported initial condition 
$\nu$. Suppose $h: [0,\infty) \to \mathbf{R}$ has two bounded continuous 
derivatives and $\phi: \mathcal{R}^{d} \to [0,\infty) $ is continuous with compact support. 
Fix $0 < t_0 < t_1$ and a bounded $\sigma( u_s(dx): 0 \leq s \leq t_0)$ variable $Z_{t_0}$. Then 
\begin{eqnarray}  
\lefteqn{  E \left[ Z_{t_0} h \left( \int \! \int \phi(x-y) u_{t_1}(dx) \nu (dy) \right) \right] - 
E \left[ Z_{t_0} h \left( \int \! \int \phi(x-y) u_{t_0}(dx) v_{t_1-t_0} (dy) \right) \right]} \nonumber\\
&=& \frac{\kappa^2}{2} E \left[ Z_{t_0} \int^{t_1}_{t_0} \int_{\mathbf{R}^{4d}}  h^{\prime \prime} 
\left(\int \! \int \phi(x-y) u_{s}(dx) v_{t_1-s}(dy)\right) \phi (x_1- y_1) \phi (x_2- y_2)
\right.  \nonumber \\
&& \hspace{.4in} \cdot \left.  \left( 
\frac{1}{|x_1-x_2|^2}- \frac{1}{|y_1-y_2|^2} \right) 
u_{s}(dx_1)v_{t_1-s}(dy_1) u_{s}(dx_2) v_{t_1-s}(dy_2) ds \right].  \label{4.2}
\end{eqnarray}
\end{Lemma}

\noindent
\textbf{Proof} We first establish some integrability, sufficient to ensure that the expectations on the right 
hand side of (\ref{4.2}) is finite. Using the independence of $\{u_t(dx)\}$ and $\{v_t(dx)\}$, 
the compact  support of $\phi$ and the bound on second moments in (\ref{1.7}), a lengthy but 
straightforward calculation, similar to that in Lemma \ref{momentcor1}, yields 
\begin{eqnarray}  
&& \hspace{-.3in} E \left[\int_{\mathbf{R}^{4d}} 
\phi( x_1 - y_1) \phi(x_2 - y_2) \left( \frac{1}{|x_1-x_2|^2} + \frac{1}{|y_1-y_2|^2} 
\right) u_{s}(dx_1) v_{t}(dy_1) u_{s}(dx_2) v_{t}(dy_2) \right]  \nonumber \\
& \leq & C(\phi, \mu, \nu, T) \left( s^{-(2-\alpha)_+/2} + t^{-(2-\alpha)_+/2} \right) \quad 
\mbox{ for all $s,t \leq T$.}  \label{4.3}
\end{eqnarray}
Furthermore, using the formula for first moments (\ref{1.6}), an easy calculation shows that 
\begin{equation}  \label{4.4}
E \left[ \int \! \int \phi(x-y) u_s(dx) v_t(dy) \right] \leq C(\phi, \mu, \nu, T) \quad 
\mbox{for all $0 \leq s,t \leq T$.}
\end{equation}
We now follow the standard method of duality, as explained in Ethier and Kurtz \cite{Ethier+Kurtz86} 
Section 4.4. Take $f \in \mathcal{C}^2_c $, apply Ito's formula using the martingale problem for $u_t(f)$ 
and then take expectations to obtain, for $s \geq t_0$, 
\begin{eqnarray*}
\lefteqn{ E[ Z_{t_0} h( u_s(f))] - E[ Z_{t_0} h( u_{t_0}(f))] } \\
&= & \int^s_{t_0} E \left[ Z_{t_0} \left( h^{\prime}( u_r(f)) u_r( \frac12 \Delta f) + \frac{\kappa^2}{2} h^{\prime \prime} 
(u_r(f)) \int \! \int \frac{f(x_1)f(x_2)}{|x_1 - x_2|^2} u_r(dx_1) u_r(dx_2) \right) \right] dr.
\end{eqnarray*}
Here Lemma \ref{momentcor1} implies that the local martingale arising from Ito's formula is a true 
martingale. Now take $\psi: \mathbf{R}^{2d} \to \mathbf{R}$, twice continuously 
differentiable and with compact support. Replace the deterministic function $f(x)$ by the random 
$\mathcal{C}^2_c$ function, independent of $\{u_t(dx)\}$, 
given by $f(x) = \int \psi(x,y) v_t(dy)$. Fubini's theorem and the integrability in 
(\ref{4.3}) and (\ref{4.4}) imply that, for $s \geq t_0$, 
\begin{eqnarray*}
\lefteqn{ E \left[  Z_{t_0} h \left( \int \! \int \psi(x,y) u_s(dx) v_{t}(dy)\right) \right]
 - E\left[ Z_{t_0}  h \left( \int \! \int  \psi(x,y) u_{t_0} (dx ) v_{t}(dy) \right) \right]} \\
&= & \int^s_{t_0} E \left[ Z_{t_0} h^{\prime} \left( \int \! \int \psi(x,y) u_r(dx) v_t(dy) \right) 
\int \! \int \frac12  \Delta^{(x)} \psi(x,y) u_r(dx) v_t(dy) \right] dr \\
& & + \frac{\kappa^2}{2} \int^s_{t_0}
 E \left[ Z_{t_0} h^{\prime \prime} \left( \int \int \psi(x,y) u_r(dx) v_t(dy) \right)  \right. \\
&& \hspace{.35in} \left.  \cdot 
\int_{\mathbf{R}^{4d}} \frac{\psi(x_1, y_1)\psi(x_2, y_2)}{|x_1 - x_2|^2} 
u_r(dx_1) u_r(dx_2) v_t(dy_1) v_t(dy_2) \right] dr.
\end{eqnarray*}
In a similar way, applying Ito's formula to $v_t(f)$, we obtain the decomposition 
\begin{eqnarray*}
\lefteqn{ E \left[ Z_{t_0} h \left( \int \! \int \psi(x,y) u_s(dx) v_{t}(dy)\right) \right] - 
E\left[ Z_{t_0} h \left( \int \! \int \psi(x,y) u_s(dx) \nu(dy) \right) \right]} \\
&= & \int^t_0 E \left[Z_{t_0} h^{\prime} \left( \int \! \int \psi(x,y) u_s(dx)
v_r(dy) \right) \int \! \int \frac12 \Delta^{(y)} \psi(x,y) u_s(dx) v_r(dy) \right] dr \\
& & + \frac{\kappa^2}{2} \int^t_0 E \left[ Z_{t_0} h^{\prime \prime} \left( \int \int \psi(x,y) 
u_s(dx) v_r(dy) \right) \right. \\
&& \hspace{.35in} \left. \cdot \int_{\mathbf{R}^{4d}} 
\frac{\psi(x_1, y_1)\psi(x_2, y_2)}{|y_1 - y_2|^2} u_s(dx_1) 
u_s(dx_2) v_r(dy_1) v_r(dy_2) \right] dr.
\end{eqnarray*}
Now defining 
\[
F(s,t)= E\left[ Z_{t_0} h \left( \int \! \int \psi(x,y) u_s(dx) v_t(dy) \right) \right] 
\]
the last two decompositions show that $s \to F(s,t)$ and $t \to F(s,t)$ are both absolutely continuous and 
gives expressions for their derivatives $\partial_1 F(s,t)$ and $\partial_2 F(s,t)$. Then applying Lemma 4.4.10 
from \cite{Ethier+Kurtz86} we obtain 
\begin{eqnarray}  
\lefteqn{ E \left[ Z_{t_0} h \left( \int \! \int \psi(x,y) u_{t_1}(dx) \nu (dy) \right) \right] 
- E \left[ Z_{t_0} h \left( \int \! \int f(x,y) u_{t_0} (dx) v_{t_1-t_0} (dy) \right) \right] } \nonumber  \\
& = & F(t_1,0) - F(t_0,t_1-t_0) \nonumber \\
& = & \int^{t_1}_{t_1-t_0} \partial_1 F(s,t_1-s) - \partial_2 F(s,t_1-s) ds  \nonumber \\
&=& \int^{t_1}_{t_0} E \left[ Z_{t_0}
 h^{\prime} \left( \int \! \int \psi(x,y) u_s(dx) v_{t_1-s}(dy) \right) \int \!  \int
 \frac12 ( \Delta^{(x)} - \Delta^{(y)}) \psi(x,y) u_s(dx) v_{t_1-s}(dy) \right] ds  \nonumber \\
& & + \frac{\kappa^2}{2} \int^{t_1}_{t_0} \int_{\mathbf{R}^{4d}} 
E \left[ Z_{t_0} h^{\prime \prime} \left(\int \! \int \psi(x,y) u_{s}(dx) 
v_{t_1-s}(dy)\right) \psi (x_1, y_1) \psi (x_2,y_2)
\right.  \nonumber \\
&& \qquad \hspace{.3in} \left. \cdot \left( \frac{1}{|x_1-x_2|^2}- \frac{1}{|y_1-y_2|^2}  \right)
u_{s}(dx_1)v_{t_1-s}(dy_1) u_{s}(dx_2) v_{t_1-s}(dy_2) \right] ds. \label{4.5}
\end{eqnarray}
Now suppose that $\phi: \mathbf{R}^{d} \to [0,\infty) $ 
is smooth and has compact support.  Choose a series of smooth, compactly support 
functions $\psi_n(x,y)$ satisfying $0 \leq \psi_n \uparrow 1$ as $n \to \infty$ and 
with $\partial_x \psi_n $,  $\partial_y \psi_n $, $\partial_{xx} \psi_n $, $\partial_{yy} \psi_n $ 
converging uniformly to zero. Apply  (\ref{4.5}) to the function 
$\psi(x,y) = \psi_n(x,y) \phi(x-y)$. Using $(\Delta^{(x)} - \Delta^{(y)})  \phi(x-y) = 0$ 
we may, using the integrability in (\ref{4.3}) and (\ref{4.4}), pass to 
the limit in (\ref{4.5})  to yield (\ref{4.2}). Finally we obtain the 
result for general  continuous $\phi$ by taking smooth approximations.  \qed

Now we take $\phi(x)$  a 
smooth, non-negative function on $\mathbf{R}^d$, supported 
on the unit ball $\{x \in \mathbf{R}^d: |x| \leq 1\}$ and 
satisfying $\int_{\mathbf{R}^d}\phi(x)dx=1$. Define an 
approximate identity by $ \phi_{\varepsilon}(x) = \varepsilon^{-d} 
\phi(x/\varepsilon)$. We may and shall suppose that $ 0 \leq \phi(x) \leq 2 G_{\varepsilon}(x)$ and 
hence that $\phi_{\varepsilon} \leq G_{\varepsilon^2}$.
We shall use this test function in Lemma \ref{doubletestfunction}) and the following lemma controls 
the right hand side of  (\ref{4.2}).
\begin{Lemma} \label{errorterm} 
Suppose $\{u_t(dx)\}$ and $\{v_t(dx)\}$ are independent solutions of (\ref{1.1}), with initial 
conditions $\mu, \nu $, where  $\nu$ compactly supported.  Then
\[
E \left[ \int^t_0 \int_{\mathbf{R}^{4d}} \phi_{\varepsilon} (x_1- y_1) \phi_{\varepsilon} 
(x_2- y_2)   \left| \frac{1}{|x_1-x_2|^2}- \frac{1}{|y_1-y_2|^2} \right|
u_{s} (dx_1) v_{t-s} (dy_1) u_{s} (dx_2) v_{t-s} (dy_2) ds \right]
\]
converges to zero as $\varepsilon \to 0$.
\end{Lemma}

\noindent
\textbf{Proof}
This lemma is a straightforward but lengthy consequence of the 
second moment bounds (\ref{1.7}).  Since it is this proof that 
requires the strict inequality $\kappa<(d-2)/2$ and also the requirement
that $\mu,\nu \in \mathcal{H}_{\beta}$ for some $\beta> \alpha$, we give some 
of the details.

The second moment bounds show that
show that the expectation in the statement of the lemma is bounded by 
\begin{eqnarray}
&& C \int^t_0 \int_{\mathbf{R}^{8d}} \phi_{\varepsilon} (x'_1- y'_1) \phi_{\varepsilon} 
(x'_2- y'_2) 
G_s(x_1-x'_1) G_s(x_2-x'_2) G_{t-s}(y_1-y'_1) G_{t-s}(y_2-y'_2) \nonumber \\
&& \hspace{.4in} \cdot 
\left( 1+ \frac{s^{\alpha}}{|x_1-x_2|^{\alpha} |x'_1 - x'_2|^{\alpha}} \right)
\left( 1+ \frac{(t-s)^{\alpha}}{|y_1-y_2|^{\alpha} |y'_1 - y'_2|^{\alpha}} \right) \nonumber \\
&& \hspace{.6in} \cdot \left| \frac{1}{|x'_1-x'_2|^2}- \frac{1}{|y'_1-y'_2|^2} \right|
 \mu (dx_1) \mu (dx_2) \nu(dy_1) \nu(dy_2) dx'_1 dx'_2 dy'_1 dy'_2 ds. \label{4.6}
\end{eqnarray}
The idea is to bound first the $dx'_1 dx'_2 dy'_1 dy'_2$ integral. 
We can split the $dx'_1 dx'_2 dy'_1 dy'_2$ integral into four terms by expanding the brackets
\[
\left( 1+ \frac{s^{\alpha}}{|x_1-x_2|^{\alpha} |x'_1 - x'_2|^{\alpha}} \right)
\left( 1+ \frac{(t-s)^{\alpha}}{|y_1-y_2|^{\alpha} |y'_1 - y'_2|^{\alpha}} \right).
\]
We shall show only how to treat the worst of these terms, namely
\begin{eqnarray}
&&  \int_{\mathbf{R}^{4d}} \phi_{\varepsilon} (x'_1- y'_1) \phi_{\varepsilon} 
(x'_2- y'_2) 
G_s(x_1-x'_1) G_s(x_2-x'_2) G_{t-s}(y_1-y'_1) G_{t-s}(y_2-y'_2) \label{4.7} \\
&& \hspace{.2in}  \cdot \left(  \frac{s^{\alpha}(t-s)^{\alpha}}{|x_1-x_2|^{\alpha} |x'_1 - x'_2|^{\alpha}
|y_1-y_2|^{\alpha} |y'_1 - y'_2|^{\alpha}} \right) \left| \frac{1}{|x'_1-x'_2|^2}- \frac{1}{|y'_1-y'_2|^2} \right|
dx'_1 dx'_2 dy'_1 dy'_2 ds.  \nonumber 
\end{eqnarray}
This is the term that requires the restriction on $\kappa$. The other three terms are similar but easier. 

We split the integral (\ref{4.7}) into two regions.  
First we consider  $x'_1,y'_1,x'_2,y'_2$ lying in the set 
$A_{\varepsilon}= \{ |x'_1-x'_2| \geq \varepsilon^{\gamma}, |y'_1-y'_2| \geq \varepsilon^{\gamma}\}$, where 
$\gamma \in (0,1)$ will be chosen later in the proof.  On this set,  since we may also suppose
$|x'_1-y'_1| \leq \varepsilon$ and $|x'_2 -y'_2| \leq \varepsilon$ by the support of 
$\phi_{\varepsilon}  $, we have, arguing using the mean value theorem,
\[
\left| \frac{1}{|x'_1-x'_2|^2}- \frac{1}{|y'_1-y'_2|^2} \right| \leq C(\gamma) \varepsilon^{1-3 \gamma}
\quad \mbox{for all $\varepsilon < 1/3$.}
\] 
This bound means  the integral (\ref{4.7}),  over the set $A_{\varepsilon}$, 
can be bounded by
\[
C(\gamma) \varepsilon^{1-3 \gamma - 2 \alpha \gamma} 
\frac{s^{\alpha}(t-s)^{\alpha}}{|x_1-x_2|^{\alpha}|y_1-y_2|^{\alpha}} 
G_{t+\varepsilon^2} (x_1-y_1)G_{t+\varepsilon^2} (x_2-y_2).
\]
We shall choose $\gamma>0 $ so that $1-3 \gamma - 2 \alpha \gamma >0$.
It is easy to show that this bound substituted into (\ref{4.6}) 
will vanish as $\varepsilon \downarrow 0$.
To estimate the integral (\ref{4.7}) over the complimentary set $A^c_{\varepsilon}$ 
we simply bound
\[
\left| \frac{1}{|x'_1-x'_2|^2}- \frac{1}{|y'_1-y'_2|^2} \right| \leq \frac{1}{|x'_1-x'_2|^2} + \frac{1}{|y'_1-y'_2|^2}
\]
and it becomes 
\begin{eqnarray*}
&& \hspace{-.2in}  \frac{s^{\alpha}(t-s)^{\alpha}}{|x_1-x_2|^{\alpha} |y_1-y_2|^{\alpha}} 
\int_{A^c_{\varepsilon}} G_s(x_1-x'_1) G_s(x_2-x'_2) G_{t-s}(y_1-y'_1) G_{t-s}(y_2-y'_2)  \\
&& \hspace{.1in} \cdot \phi_{\varepsilon} (x'_1- y'_1) \phi_{\varepsilon} 
(x'_2- y'_2)   |x'_1 - x'_2|^{-\alpha} |y'_1 - y'_2|^{-\alpha} \left( |x'_1-x'_2|^{-2} + 
|y'_1-y'_2|^{-2} \right) dx'_1 dx'_2 dy'_1 dy'_2 \\
 & \leq &   C  \frac{s^{\alpha}(t-s)^{\alpha}}{|x_1-x_2|^{\alpha} |y_1-y_2|^{\alpha}} 
\int_{A^c_{\varepsilon}} G_s(x_1-x'_1) G_s(x_2-x'_2) G_{t-s}(y_1-y'_1) G_{t-s}(y_2-y'_2)  \\
&& \hspace{.2in} \cdot \phi_{\varepsilon} (x'_1- y'_1) \phi_{\varepsilon} 
(x'_2- y'_2) \left( |x'_1-x'_2|^{-(2+2 \alpha)} + 
|y'_1-y'_2|^{-(2+2 \alpha)} \right) dx'_1 dx'_2 dy'_1 dy'_2.
\end{eqnarray*}
We show how to treat just the integral with the term $|x'_1-x'_2|^{-(2+2\alpha)}$, the term 
$|y'_1-y'_2|^{-(2+2\alpha)}$ being entirely similar. Note that the restriction $\kappa < (d-2)/2$ is 
simply to ensure that $2 + 2 \alpha < d$ and so the pole $|z|^{-(2+2 \alpha)}$ is integrable on $\mathbf{R}^d$. We may choose $\delta \in (2+2 \alpha, (2+\alpha + \beta) \wedge d)$. Then, using the bound $\phi_{\varepsilon} \leq 2 G_{\varepsilon^2}$, we can do the $dy'_1 dy'_2$ integrals to see that 
\begin{eqnarray*}
&& \hspace{-.2in} \int_{A^c_{\varepsilon}} G_s(x_1-x'_1) G_s(x_2-x'_2)
 G_{t-s}(y_1-y'_1) G_{t-s}(y_2-y'_2) \\
&& \hspace{.3in} \cdot \phi_{\varepsilon} (x'_1- y'_1) \phi_{\varepsilon} 
(x'_2- y'_2)  |x'_1-x'_2|^{-(2+2 \alpha)} dx'_1 dx'_2 dy'_1 dy'_2 \\
& \leq  &  \int_{\{|x'_1-x'_2| \leq \varepsilon\}} G_s(x_1-x'_1) G_s(x_2-x'_2) \\
&& \hspace{.3in} \cdot G_{t-s+\varepsilon^2}(y_1-x'_1) 
G_{t-s+\varepsilon^2}(y_2-x'_2)   |x'_1-x'_2|^{-(2+2 \alpha)} dx'_1 dx'_2  \\
& \leq & C \varepsilon^{(\delta-2-2\alpha)\gamma} 
\int_{\mathbf{R}^{2d}} G_s(x_1-x'_1) G_s(x_2-x'_2)  \\
&& \hspace{.3in} \cdot  G_{t-s+\varepsilon^2}(y_1-x'_1)
 G_{t-s+\varepsilon^2}(y_2-x'_2)   |x'_1-x'_2|^{-\delta} dx'_1 dx'_2 .
\end{eqnarray*}
We now split into two cases: $s \leq t/2$  and $s \geq t/2$. When $s \leq t/2$ we have the  bound 
\[
G_{t-s+\varepsilon^2}(y_1 - x'_1) G_{t-s+\varepsilon^2}(y_2 - x'_2) \leq C(a, \nu, t) \exp( - a|x'_1| - a|x'_2|),
\quad \mbox{for $y_1, y_2 \in $supp($\nu$).}
\]
So, when $ s \leq t/2$,  
\begin{eqnarray*}
\lefteqn{\int_{\mathbf{R}^{2d}} G_s(x_1-x'_1) G_s(x_2-x'_2)  G_{t-s+\varepsilon}(y_1-x'_1) G_{t-s+\varepsilon}(y_2-x'_2)   |x'_1-x'_2|^{-\delta} dx'_1 dx'_2 } \\
& \leq & C(a, \nu, t) \int_{\mathbf{R}^{2d}} G_s(x_1-x'_1) G_s(x_2-x'_2) 
\exp( - a|x'_1| -a|x'_2|)  |x'_1-x'_2|^{-\delta} dx'_1 dx'_2 \\
& \leq & C(a,\nu,t) \exp(-a|x_1| -a |x_2| ) \left( |x_1-x_2|^{-\delta} \wedge s^{-\delta/2} \right) \\
& \leq & C(a,\nu,t) \exp(-a|x_1| -a |x_2| ) |x_1-x_2|^{-\beta + \alpha} s^{-(\delta-\beta+\alpha)/2},
\end{eqnarray*}
using the tricks from Lemma \ref{momentcor1} for this last inequality. 
Combining all these bounds one has, when substituting the integral (\ref{4.7}) over the region
$A^c_{\varepsilon}$ into (\ref{4.6}), and considering only  the time interval $[0,t/2]$, the estimate
\[
C(a, \nu, t) \varepsilon^{(\delta-2-2\alpha)\gamma}  \int^{t/2}_0 \int_{\mathbf{R}^{4d}}  \left(  
\frac{ s^{-(\delta-\beta-\alpha)/2} }{|x_1-x_2|^{\beta} |y_1-y_2|^{\alpha}} \right) e^{-a|x_1| -a |x_2|} 
\mu(dx_1) \mu(dx_2) \nu(dy_1) \nu(dy_2) ds.
\]
Choosing $a$ so that $e^{-a|x|} \mu(dx) \in \mathcal{H}^0_{\beta}$,
the integral is finite and so this expression vanishes as $\varepsilon \downarrow 0$. 
The integral over $[t/2,t]$ is treated in a fairly similar way
using the assumption that $\nu \in \mathcal{H}_{\beta}$.  \qed

To deduce Proposition \ref{dualityrelation} from Lemmas \ref{doubletestfunction} and \ref{errorterm} is 
easy. By a simple approximation argument it is enough to prove (\ref{4.1}) for $h$ with two 
bounded continuous derivatives. We apply the approximate duality relation (\ref{4.2}),
using $0=t_0 < t_1=t$ and $Z_{t_0}=1$, to the  function $\phi_{\varepsilon}$.  
Then take $\varepsilon \to 0$ and use the control on the error term in  Lemma \ref{errorterm} to 
obtain the result.  

We show two consequences of the duality relation and its proof. 
\begin{Corollary} \label{dualitycor1}
Solutions to (\ref{1.1}) are unique in law and we let $Q_{\mu}$
denote the law of solutions started at $\mu \in \mathcal{H}_{\alpha+}$.
\end{Corollary}
\noindent
\textbf{Proof} First suppose that $\{u_t(dx)\}$ and $\{v_t(dx)\}$ are two solutions with 
the same deterministic initial condition $\mu$. 
Construct a third solution $\{w_t(dx)\}$,
independent of $\{u_t(dx)\}$ and $\{v_t(dx)\}$ and with initial condition $w_0(dx) = f(x) dx$ 
for some non-negative, continuous, compactly supported function $f$. Then apply the 
approximate duality relation (\ref{4.2}), with $0=t_0 < t_1 =t$ and $Z_{t_0}=1$,
to the pair $\{u_t(dx)\}$ and $\{w_t(dx)\}$ and to the pair $\{v_t(dx)\}$ and $\{w_t(dx)\}$, 
using the function $\phi_{\varepsilon}$.
Subtracting the two approximate duality relations we see that 
\[
E\left[ h\left( \int \! \int \phi_{\varepsilon} (x-y) u_t(dx) f(y) dy
\right) \right] - E\left[ h\left( \int \! \int \phi_{\varepsilon} (x-y)
v_t(dx) f(y) dy \right) \right] 
\]
equals the sum of two error terms, both of which converge to zero as
 $\varepsilon \to 0$ by Lemma \ref{errorterm}. Hence $E[ h(u_t(f))] = E[ h(
v_t(f)) ] $ for all such $f$ and for all suitable $h$. Choosing $h(z) =
\exp(-\lambda z)$ we obtain equality of the Laplace functionals of $u_t(dx)$ and 
$v_t(dx)$ and hence equality of the one dimensional distributions. 

Now we use an induction argument to show that the finite dimensional
distributions agree. Suppose the $n$-dimensional distributions have been shown to agree.
Choose $0 \leq s_1 < s_2 \ldots < s_{n+1}$ and set $t_1=s_{n+1}, t_0 = s_n$. 
Then apply the approximate duality relation (\ref{4.2})  to the pair 
$\{u_t(dx)\}$ and  $\{w_t(dx)\}$ with $Z_{t_0} = \prod_{i=1}^n \exp( - u_{s_i}(f_i))$ 
for compactly supported $f_i \geq 0$. Also apply the approximate duality relation 
(\ref{4.2})  to the pair  $\{v_t(dx)\}$ and  $\{w_t(dx)\}$ with 
$Z_{t_0} = \prod_{i=1}^n \exp( - v_{s_i}(f_i))$.
Subtracting the two approximate duality relations, use the equality of the
$n$-dimensional distributions and let $\varepsilon \downarrow 0$ to obtain
equality of the $n+1$-dimensional distributions, completing the induction.
Since the paths have continuous paths the finite dimensional distributions 
determine the law.

For general initial conditions $u_0$ we let $P_{\mu}$ be a regular conditional
probability given that $u_0 = \mu$. It is not difficult to check that for
almost all $\mu$ (with respect to the law of $u_0$) the process $\{u_t(dx)\}$
is a solution to (\ref{1.1}) started at $\mu$ under $P_{\mu}$. (The 
moment conditions carry over under the regular conditional probability and
these allow one to reduce to a countable family of test functions in the 
martingale problem). By the argument above the law of $\{u_t(dx)\}$ under
the conditional probability 
$P_{\mu}$ is uniquely determined (for almost all $\mu$). This in turn determines  
the law of $\{u_t(dx)\}$.  \qed

\begin{Corollary} \label{dualitycor2} 
For any bounded Borel measurable $H: C([0,\infty),\mathcal{M}) \to  \mathbf{R}$ the map 
$\mu \rightarrow Q_{\mu}[H]$, the integral of $H$ with respect to $Q_{\mu}$, is measurable from $\mathcal{H}_{\alpha+} $ to $\mathbf{R}$.

The set of laws $\{Q_{\mu}: \mu \in \mathcal{H}_{\alpha+}\}$ forms a Markov family, 
in that for any solution $\{u_t(dx)\}$ to (\ref{1.1}), for any bounded measurable
$H: C([0,\infty),\mathcal{M}) \to \mathbf{R}$, and for any $t \geq 0$
\[
E\left[ H( u_{t +\cdot}) | \mathcal{F}_{t}\right] = Q_{u_{t}}
\left[ H \right], \; \; \mbox{almost surely.}
\] 
\end{Corollary}
\noindent
\textbf{Proof.} We use the methods of  Theorem 4.4.2 of  Ethier and Kurtz \cite{Ethier+Kurtz86}. We were unable to directly apply these results, but with a little adjustment
the methods apply to our case and we point out the key changes needed. 

We only allow initial conditions in the strict subset $\mathcal{H}_{\alpha+}$ of all Radon measures, 
and do not yet know that the process takes values in this subset. But by
restricting to the ordinary Markov property it is enough to know that
$P(u_t(dx) \in \mathcal{H}_{\alpha+})=1$ for each fixed $t$, and this follows from
Lemma \ref{momentcor1} part ii). 

The measurability of $\mu \to Q_{\mu}[H]$ can often be established for
martingale problems by establishing it as the inverse of a suitable Borel
bijection (see \cite{Ethier+Kurtz86} Theorem 4.4.6). We do not use this method as 
$\mathcal{H}_{\alpha+}$ is not complete under the vague topology. However the
measurability can be established directly as follows. It is enough, by a
monotone class argument, to consider $H$ of the form 
$H(\omega) = \prod_{i=1}^n h_i(\omega_{t_i}(f_i))$
for bounded continuous functions $h_i$, for $f_i \in \mathcal{C}_c$, for 
$0 \leq t_1<t_2 \ldots t_n$ and for $n \geq 1$. But for such $H$ we can write,
using the construction of solutions from Section \ref{s4}, 
\[
Q_{\mu}[H] =  E \left[ \prod_{i=1}^n h_i \left( \sum_{n=0}^{\infty}
I^{(n)}_{t_i} (f_i, \mu) \right) \right] 
=  \lim_{N \to \infty} E \left[ \prod_{i=1}^n h_i \left( \sum_{n=0}^{N}
I^{(n)}_{t_i} (f_i, \mu) \right) \right] .
\]
For each $N < \infty$ the integrands 
$ \sum_{n=0}^{N} I^{(n)}_{t_i} (f_i, \mu) $ are, by the definition of the
maps $I^{(n)}(f,\mu)$, continuous in $\mu$. So $Q_{\mu}[H]$ is the limit
of continuous maps on $\mathcal{H}_{\alpha+}$. 

We can now follow  the method of in Theorem 4.4.2 part c) in 
Ethier and Kurtz \cite{Ethier+Kurtz86} in  the proof of the Markov property.
The only important change in the argument from
Ethier and Kurtz is that we have uniqueness in law for solutions to (\ref{1.1}), 
and this requires the moment bounds  (\ref{1.6}) and (\ref{1.7})
to hold as well the martingale problem (\ref{1.4})  and 
(\ref{1.5}). The key point is to show that, for any $t>0$, 
the process $\{u_{t+\cdot}(dx)\}$ satisfies these moment bounds. For this 
 it is enough to show for all $f: \mathbf{R}^d \to [0,\infty)$ and $0 < s < t$
\[
E \left[ u_t(f) | \sigma(u_r(dx): r \leq s)  \right] =  \int \int  G_{t-s}(x-x') f(x') u_s(dx),
\]
and there exists $C$, depending only on the dimension $d$ and $\kappa$, so that 
\begin{eqnarray*}  
\lefteqn{ E  \left[ \left. \left(  \int f(x) u_t(dx) \right)^2 \right| \sigma( u_r(dx): r \leq s) \right]} \\
& \leq & C  \int_{\mathbf{R}^{4d}} G_{t-s} (x-x') G_{t-s}(y-y') f(x') f(y') 
\left(1+  \frac{(t-s)^{\alpha} }{|x-y|^{\alpha}|x'-y'|^{\alpha}} \right) 
u_s(dx) u_s(dy) dx' dy'.  
\end{eqnarray*}
By uniqueness in law it is enough to prove these bounds for the solutions
constructed via chaos expansions in Section \ref{s4}. 
It is also enough to prove these bounds for  $f \in \mathcal{C}_c$.
The first moment follows 
from the fact that $E[ I^{(n)}_t(f,\mu) | \mathcal{F}_s] = I^{(n)}_s(G_{t-s}f, \mu)$
and the convergence of the series (\ref{3.2}).
For the second moment bound
we use the approximations $u^{(\varepsilon)}_t$ introduced in Section \ref{s4},
for which we know $u^{(\varepsilon)}(f) \to u_t(f)$ in $L^2$.
Fix $0 < s_1 < \ldots < s_n \leq s$, $f_1, \ldots , f_n \in \mathcal{C}_c$ and a bounded
continuous function $h: \mathbf{R}^n \to \mathbf{R}$. Then, using the Markov
property of the approximations $u^{(\varepsilon)}_t$,
\begin{eqnarray}
\lefteqn{E \left[ (u_t(f))^2 h(u_{s_1}(f_1), \ldots, u_{s_n}(f_n)) \right]} \nonumber \\
& = & \lim_{\varepsilon \downarrow 0} 
E \left[ (u^{(\varepsilon)}_t(f))^2 
h(u_{s_1}(f_1), \ldots,  u_{s_n}(f_n)) \right] \nonumber  \\
& \leq & C \lim_{\varepsilon \downarrow 0}  
E \left[\int_{\mathbf{R}^{4d}} G_{t-s} (x-x') G_{t-s}(y-y') f(x') f(y') \right.  \nonumber  \\
&& \hspace{.3in} \cdot
\left. \left(1+  \frac{(t-s)^{\alpha} }{|x-y|^{\alpha}|x'-y'|^{\alpha}} \right) 
u^{(\varepsilon)}_s(dx) u^{(\varepsilon)}_s(dy) dx' dy' \;
h(u_{s_1}(f_1), \ldots,  u_{s_n}(f_n)) \right] \nonumber  \\
& = & C 
E \left[\int_{\mathbf{R}^{4d}} G_{t-s} (x-x') G_{t-s}(y-y') f(x') f(y') \right. \label{4.8} \\ 
&& \hspace{.5in} \cdot 
\left. \left(1+  \frac{(t-s)^{\alpha} }{|x-y|^{\alpha}|x'-y'|^{\alpha}} \right) 
u_s(dx) u_s(dy) dx' dy' \;  h(u_{s_1}(f_1), \ldots,  u_{s_n}(f_n)) \right]. \nonumber 
\end{eqnarray}
The last equality follows by the convergence $u^{(\varepsilon)}(f) \to u_t(f)$ for
compactly supported $f$ and an approximation argument using the uniform second moment
bounds on $u^{(\varepsilon)}_s$ and $u_s$. 
The inequality (\ref{4.8}) implies the desired second moment bound and completes the proof. 
\qed
%
%
%
\section{Death of solutions}
%
%
\label{s6} \setcounter{equation}{0}

We adapt a method from the particle systems literature to 
study questions of extinction.  Liggett and Spitzer used this technique, described in 
Chapter IX, Section 4 of \cite{Liggett85}, to study analogous questions for linear 
particle systems. The corresponding result  for certain linear particle systems, indexed on 
$\mathbf{Z}^d$ and with noise that is white in space, is that death of solutions 
occurs in dimensions $d=1,2$ for all $\kappa$, and in dimensions $d \ge 3$ for 
sufficiently large $\kappa$. The long range correlations of our noise lead to different 
behavior, an increased chance of death, and death occurs for all the values of
$d \geq 3$ and $\kappa$ that we are considering. However our basic estimate in the proof
of Proposition \ref{liggett} below leaves open the possibility that the death is extremely slow.

We start by considering initial conditions with finite total mass. To study the evolution of 
the total mass we want to use the test function $f=1$ in the martingale problem. The next lemma 
shows this is possible by approximating $f$ by suitable compact support test functions.

\begin{Lemma} \label{totalmass} 
Suppose that $\{u_t(dx)\}$ is a solution 
to (\ref{1.1}) started at $\mu \in \mathcal{H}^0_{\alpha}$. Then the total mass 
$\{u_t(1): t \geq 0 \}$ is a continuous martingale with 
\[ 
\langle u(1) \rangle_t = \int^t_0 \int \int  \frac{u_s(dx) \, u_s(dy)}{|x-y|^2} ds. 
\]
\end{Lemma}

\noindent
\textbf{Proof.} We first check that the assumptions on the initial condition imply that 
$E[u_t(1)^2] < \infty$. The bound on second moments  (\ref{1.7}) implies that 
\begin{eqnarray} 
E \left[ u_t(1)^2 \right]
& \leq & C \int_{\mathbf{R}^{4d}} G_t(x-x') G_t(y-y')   
\left(1+ \frac{t^{\alpha}}{ |x-y|^{\alpha}|x'-y'|^{\alpha}} \right) 
dx' dy' \mu(dx) \mu(dy) \nonumber \\
& \leq & C \int_{\mathbf{R}^{2d}} \left(1+ \frac{t^{\alpha/2}}{ |x-y|^{\alpha}}\right) 
\mu(dx) \mu(dy)   \quad \mbox{using (\ref{1.9})} \nonumber \\
& \leq & C (1 + t^{\alpha/2}) \|\mu \|_{\alpha}^2. \label{5.1}
\end{eqnarray}
We may find $f_n \in C_c^2 (\mathbf{R}^d)$ so that $0 \leq f_n \uparrow 1$ and 
$\| \Delta f_n\|_{\infty} \downarrow 0$ as $n \rightarrow \infty$. Applying Doob's 
inequality we have, for any $T \geq 0$, 
\begin{eqnarray*}
\lefteqn{ E \left[ \sup_{t \leq T} | z_t(f_n)-z_t(f_m) |^2 \right]} \\
& \leq & C E \left[ | z_T(f_n)-z_T(f_m) |^2 \right] \\
& = & C E\left[ \left|u_T(f_n-f_m) - \mu(f_n-f_m) - \frac12 \int^T_0 u_s 
(\Delta f_n - \Delta f_m) ds \right|^2 \right] \\
& \leq & C \left( E \left[ u_T(f_n-f_m)^2 \right] + \mu (f_n-f_m)^2 + \frac12 
( \| \Delta f_n\|_{\infty} + \| \Delta f_m\|_{\infty}) 
E \left[ \left(\int^T_0 u_s(1)ds \right)^2 \right] \right) .
\end{eqnarray*}
This expression is seen to converge to zero as $n,m \rightarrow \infty$ by using dominated 
convergence and the bound in (\ref{5.1}). From this we can deduce that, along a 
subsequence, $z_t(f_n)$ converges uniformly on compacts to a continuous martingale. Also 
\[
E \left[ \sup_{t \leq T} \left| \int^t_0 u_s(\Delta f_n) ds \right|^2 \right]  
\leq \|\Delta f_n \|^2_{\infty} T \, \int^T_0 E[u^2_s(1)]ds \rightarrow 0. 
\]
Since 
\[
z_t(f_n) + \int^t_0 u_s( \frac12 \Delta f_n) = u_t(f_n) - \mu(f_n) \rightarrow u_t(1)-\mu(1) 
\]
we can conclude that $u_t(1)$ is a continuous martingale. Moreover we claim that 
\begin{eqnarray}  \label{5.2}
\lefteqn{ E \left[ \sup_{t \leq T} \left| \int^t_0 \int \int  
\frac{u_s(dx)\, u_s(dy)}{|x-y|^2} ds - \langle z(f_n) \rangle_t \right| \right]} \\
& = & E \left[ \int^T_0 \int \int \frac{1-f_n(x)f_m(y)}{|x-y|^{2}} u_s(dx)
u_s(dy) ds \right] \rightarrow 0 \quad \mbox{as $n \rightarrow \infty$.} 
\nonumber
\end{eqnarray}
This follows by dominated convergence and the bound 
\begin{eqnarray*}
E \left[\int^T_0 \int \int \frac{u_s(dx) \, u_s(dy)}{|x-y|^2} ds \right] 
&=& \lim_{n \rightarrow \infty} E \langle z(f_n) \rangle_T   \\
& = & \lim_{n \rightarrow \infty} E \left[ \left( u_T(f_n)- \mu(f_n) - \int^t_0 
u_s(\Delta f_n) ds \right)^2 \right]   \\
& = & E [ u_T(1)^2] - \mu(1)^2 < \infty.  
\end{eqnarray*}
Using (\ref{5.2}) it is now straightforward to identify the quadratic variation of $u_t(1)$ as in the statement 
of the lemma.  \qed
\begin{Proposition} \label{liggett} 
Suppose that $\{u_t(dx)\}$ is a solution 
to (\ref{1.1}) started at $\mu \in \mathcal{H}_{\alpha+} \cap \mathcal{H}^0_{\alpha}$. Then 
\[
\lim_{t \to \infty}u_t(1)=0 \quad \mbox{almost surely.} 
\]
\end{Proposition}

\noindent
\textbf{Proof.} The previous lemma shows that the process $u_t(1)$ is a non-negative 
martingale and hence converges almost surely. We will show that 
\begin{equation}  \label{5.3}
\lim_{t\to\infty}E\left[u_t(1)^{1/2}\right]=0
\end{equation}
which then implies that the limit of $u_t(1)$ must be zero. We consider first the case that 
$\mu$ is compactly supported inside the ball  $B(0,K)$. We let $C_t=B(0,R_t)$ be the closed 
ball with radius 
\[
R_t=K + \left[c_8(t \vee 1)\log\log (t \vee 4)\right]^{1/2} 
\]
where $c_8$ is a fixed constant satisfying $c_8>4$. We write $C_t^c$ for the complement of 
this ball. Let $\tau_0$ be the first time $t \geq 0$ that $u_t(1)=0$.  (In a later section we shall 
show that $P(\tau_0=\infty)=1$ whenever $\mu(1)>0$ but we do not need to assume this here.) 
Using  Ito's formula, and labeling any local martingale terms by $dM$, we find that for $t<\tau_0$, 
\begin{eqnarray}  \label{5.4}
du_t(1)^{1/2}&=&dM_t-\frac{\kappa^2}{8} u_t(1)^{1/2} \int \int 
\frac{u_t(dx)u_t(dy)}{u_t(1)^2|x-y|^2} \; dt \\
& \leq & dM_t- \frac{ \kappa^2}{32R_t^2} u_t(1)^{1/2} \int_{C_t} \int_{C_t} 
\frac{u_t(dx)u_t(dy)}{u_t(1)^2} \; dt  \nonumber \\
& = & dM_t-\frac{ \kappa^2}{32R_t^2} u_t(1)^{1/2} \left(1-\int_{C_t^c} 
\frac{u_t(dx)}{u_t(1)}\right)^2 dt  \nonumber \\
& \leq & dM_t-\frac{ \kappa^2}{32R_t^2} u_t(1)^{1/2} \left(1-2\int_{C_t^c} 
\frac{u_t(dx)}{u_t(1)}\right) dt  \nonumber \\
& \leq & dM_t-\frac{ \kappa^2}{32R_t^2} u_t(1)^{1/2}
\left(1-2\left[\int_{C_t^c} \frac{u_t(dx)}{u_t(1)}\right]^{1/2}\right) dt 
\nonumber \\
& = & dM_t- \frac{ \kappa^2}{32R_t^2} u_t(1)^{1/2} dt +
 \frac{ \kappa^2}{16R_t^2} u_t(C_t^c)^{1/2} dt.  \nonumber
\end{eqnarray}
The local martingale term in (\ref{5.4}) is given by
$dM_t = (1/2) u_t(1)^{-1/2}du_t(1)$ and is reduced by the stopping times 
$\tau_{1/n} = \inf \{t: u_t(1) \leq 1/n\}$. So applying (\ref{5.4}) at the time 
$t \wedge \tau_{1/n}$ and taking expectations we obtain 
\[
E[ u_{t \wedge \tau_{1/n}} (1)^{1/2}] \leq  \mu(1)^{1/2} 
- E \left[ \int^{t \wedge \tau_{1/n}}_0 \frac{\kappa^2}{32 R_s^2} u_s(1)^{1/2} ds \right] 
+ E \left[ \int^{t \wedge \tau_{1/n}}_0 \frac{\kappa^2}{16 R_s^2} u_s(C_s^c)^{1/2} ds \right].
\]
Letting $n \rightarrow \infty$, using monotone convergence and the  moments established in 
(\ref{5.1}), we obtain the same inequality with $\tau_{1/n}$ replaced by $\tau_0$. Since 
the paths of a non-negative local martingale must remain at zero after hitting zero  we may further 
replace $t \wedge \tau_0$ by $t$ in the inequality. Defining $\eta_t=E \left[u_t(1)^{1/2}\right] $ we 
therefore have 
\begin{equation}  \label{5.5}
\eta_t \leq \eta_0 -\int_{0}^{t} \frac{\kappa^2 }{32R_s^2} \eta_s ds +
\int^t_0 \frac{\kappa^2}{16R_s^2} E\left[u_s(C_s^c)^{1/2}\right] ds
\end{equation}
The aim is to estimate the expectation in this inequality and to show that
it implies that $\eta_t \rightarrow 0$. Let 
\[
\Xi(s,t)=\exp\left(-\int_{s}^{t}\frac{\kappa^2}{32R_r^2}dr \right). 
\]
It follows from the definition of $R_t$ that $\int_{s}^{\infty} 
\frac{\kappa^2}{32R_r^2}dr =\infty $ for any $s \geq 0$ and so, 
for any $s \ge 0$, 
\begin{equation}  \label{5.6}
\lim_{t\to\infty}\Xi(s,t)=0.
\end{equation}
Applying Gronwall's inequality to (\ref{5.5}), we obtain 
\begin{equation}  \label{5.7}
\eta_t \leq \Xi(0,t)\eta_0+\int_{0}^{t}\Xi(s,t) \frac{ \kappa^2}{16R_s^2}
E\left[u_s(C_s^c)^{1/2}\right] ds.
\end{equation}
If we show that 
$\int_{0}^{\infty} \frac{\kappa^2}{16R_s^2} E\left[u_s(C_s^c)^{1/2}\right] ds< \infty$
 it then follows that $\eta_t \downarrow 0$ as $t \rightarrow \infty$ (use $0 \leq \Xi(s,t) \leq 1$, 
(\ref{5.6}) and dominated convergence). Using the Cauchy-Schwartz 
inequality and the formula for first moments, we obtain 
\begin{eqnarray}  
\left(E\left[u_t(C_t^c)^{1/2}\right]\right)^2  & \leq & E\left[u_t(C_t^c)\right] \nonumber  \\
& \leq & \int_{C_t^c}\int_{\mathbf{R}^d} G_t(x-y) \mu (dx)dy  \nonumber \\
& \leq & C \mu(1)\int_{R_t-K}^{\infty}(2\pi t)^{-d/2} \exp (-r^2/2t) r^{d-1}dr  \nonumber \\
&=& C \mu(1)\int_{(R_t-K)/\sqrt t}^{\infty} \exp (-s^2/2)s^{d-1}ds  \nonumber \\
& \leq & C \mu(1) \exp (- (R_t-K)^2 /2t) \left( 1+ ( R_t -K / \sqrt{t})^{d-1} \right)  \nonumber \\
& \leq & C \mu(1)\left[\log (t \vee 4)\right]^{-c_8/2} \left[ \log \log (t
\vee 4) \right]^{(d-1)/2}. \label{5.8} 
\end{eqnarray}
Here we have used the following standard inequality: by the change of
variables $y=x+z$ we find 
\begin{eqnarray*}
\int_{x}^{\infty}\exp (-y^2/2) y^{d-1}dy & \leq &
C \exp (-x^2/2 ) \int_{0}^{\infty}\exp(-z^2/2) (x^{d-1}+z^{d-1})dz \\
& \leq & C(1 +x^{d-1}) \exp (-x^2/2), \quad \mbox{when $x \geq 0$.}
\end{eqnarray*}
Finally we use (\ref{5.8}) to derive the following: 
\[
\int_{0}^{\infty}\frac{\kappa^2}{16R_s^2} E\left[u_s (C^c_s)\right]^{1/2}ds
\leq  C(\kappa, \mu) \left(1 + \int^{\infty}_4 
\frac{(\log \log s)^{(d-5)/4}}{s(\log s)^{c_8/4}} ds \right) < \infty.
\]
This completes the proof in the case $\mu$ is compactly supported.  In the
general case we fix $\varepsilon >0$ and split the initial condition so that 
$\mu = \mu^{(1)} + \mu^{(2)}$ where $\mu^{(1)}(1) \leq \varepsilon $ 
and $\mu^{(2)}$ is compactly supported. By uniqueness in law we 
may consider any
solution with initial condition $\mu$ and we choose to construct one as
follows: let $u^{(1)}, u^{(2)}$ be the strong solutions, as constructed in
Section 3, with respect to the same noise and with initial conditions 
$\mu^{(1)}, \mu^{(2)}$ and set $u=u^{(1)}+u^{(2)}$. It is easy to check  
that $u$ is a solution starting at $\mu$, which is a statement of the linearity
of the equation. Using Cauchy-Schwartz and the formula for first moments 
(\ref{1.6}) we have 
\begin{eqnarray*}
E \left[ u_t(1)^{1/2}  \right] & = & E \left[ \left( u^{(1)}_t(1) +
u^{(2)}_t(1)\right)^{1/2} \right] \\
& \leq & E \left[ u^{(1)}_t(1)^{1/2} \right] +E \left[ u^{(2)}_t(1)^{1/2}
\right] \\
& \leq & E \left[ u^{(1)}_t(1)\right]^{1/2} +E \left[ u^{(2)}_t(1)^{1/2}
\right] \\
& = & \varepsilon^{1/2} +E \left[ u^{(2)}_t(1)^{1/2} \right].
\end{eqnarray*}
Thus (\ref{5.3}) follows from the compactly supported case and the proposition 
is proved.  \qed

\textbf{Proof of Theorem \ref{t4}.}  Firstly the case of an
initial condition with finite total mass. If 
$P(u_0 \in \mathcal{H}^0_{\alpha+})=1$ then, 
\[
P(u_t(1) \rightarrow 0) = \int_{\mathcal{H}^0_{\alpha+}} Q_{\mu} ( U_t(1)
\rightarrow 0) P(u_0 \in d\mu) =1. 
\]
Secondly the case where of an initial condition that has locally bounded
intensity. For such $u_0$ we have, using the first moment formula, 
\begin{equation}  \label{5.9}
E[u_1(\phi)] = E \left( \int \int G_1 (x-z) \phi(x) u_0(dz) dx \right) \leq
C \int \phi(x) dx
\end{equation}
for some constant $C < \infty$. That is $E[u_1(dx)] \leq C L(dx)$ where we
write $L(dx)$ for Lebesgue measure.
Fix a bounded set $A$. By the linearity of the equation the map $\mu
\rightarrow Q_{\mu}(U_t(A) \wedge1) $ is increasing in $\mu$. Moreover it is
concave in $\mu$. Indeed if $u^{\mu}_t(dx)$ and $u^{\nu}_t(dx)$ are 
solutions started from $\mu$ and $\nu$, with respect to the same noise,
then, by linearity and the concavity of $f(z) = z \wedge 1$, 
\begin{eqnarray*}
Q_{\theta \mu + (1-\theta) \nu} (U_t(A) \wedge 1) & = & E \left[ (\theta
u^{\mu}_t(A) + (1-\theta) u^{\nu}_t(A)) \wedge 1 \right] \\
& \geq & E \left[ \theta (u^{\mu}_t(A) \wedge 1) + (1-\theta) (u^{\nu}_t(A) \wedge 1) \right] \\
& = & \theta Q_{\mu} (U_t(A) \wedge 1) + (1-\theta) Q_{\nu} (U_t(A) \wedge1).
\end{eqnarray*}
Then 
\begin{eqnarray*}
E [u_{t+1} (A) \wedge 1] & = & \int_{\mathcal{H}_{\alpha+}} Q_{\mu} (U_t(A)
\wedge 1) P( u_1 \in d \mu) \quad \mbox{(by the Markov property)} \\
& \leq & Q_{E[u_1(dx)]} (U_t(A) \wedge 1) \quad \mbox{(by Jensen's
inequality)} \\
& \leq & Q_{C L(dx)} (U_t(A) \wedge 1) \quad \mbox{(using (\ref{5.9}))} \\
& = & Q_{I_A} (C U_t(1) \wedge 1) \quad \mbox{(by self duality)}
\end{eqnarray*}
which converges to zero by Proposition \ref{liggett}. This completes the
proof of Theorem \ref{t4}. \qed
%
%
\section{Support Properties} 
%
%
\label{s7} \setcounter{equation}{0} 

In this section we establish the various properties  listed in  Theorem \ref{t5}. 
%
%
\subsection{Dimension of Support} \label{s7.1}
%
%
We can apply Frostman's Lemma (see \cite{Falconer85} Corollary 6.6) to obtain a
lower bound on the Hausdorff dimension of supporting sets for solutions
$u_t(dx)$. Indeed Lemma \ref{momentcor1} part ii) and Frostman's Lemma imply that
any non-empty Borel supporting set for the measure $u_t(dx) $, at a fixed 
$t>0$,  must, almost surely, have dimension at least $d-\alpha $. We prove in Subsection
\ref{s7.2} that if $\mu  \neq 0$ then $u_t \neq 0$ almost surely. 
This establishes the fixed $t$ result in Theorem \ref{t5} i). 
We now show a weaker lower bound that holds at all times.

\begin{Proposition} \label{dimension} 
Suppose that $\{u_t(dx)\} $ is a solution to (\ref{1.1}).
\begin{description}
\item[i.)] If $ P(u_0 \in \mathcal{H}_{(d-2-\alpha)-})=1 $ then $P(u_t \in \mathcal{H}_{(d-2-\alpha)-} \; 
\mbox{for all $t \geq 0$}) = 1$. Indeed, for some $a$, 
\[
P \left(\mbox{There exists $a$ so that } \sup_{s \leq t} \| \mu e^{-a|x|} \|_{\beta} 
< \infty \right) =1 \quad \mbox{ for all $t \geq 0$ and $\beta < d-2-\alpha$.}
\]
\item[ii.)] For any initial condition we have $P(u_t \in \mathcal{H}_{(d-2-\alpha)-} \; 
\mbox{for all $t > 0$}) = 1$. 
\end{description}
\end{Proposition}

\noindent
\textbf{Remarks}

\textbf{1.} Since $\mathcal{H}_{(d-2-\alpha)-} \subseteq \mathcal{H}_{\alpha+} $ 
(which requires $ \kappa < (d-2)/2)$) we also have, for any initial condition,
$P(u_t \in \mathcal{H}_{\alpha+} \; \mbox{for all $t > 0$}) = 1$. 

\textbf{2.} Using Frostman's Lemma, part ii) of this proposition implies that, at all times $t>0$, a 
Borel set $A_t$ that supports $u_t(dx)$ must have Hausdorff dimension at least $d-\alpha-2$. 

\textbf{3.} The idea behind the proof of Proposition \ref{dimension} is to show, for suitable values 
of $p$, that the process $ S^{(\rho)}_t = \int \int u_t(dx)u_t(dy) / |x-y|^{\rho} $ is a non-negative 
supermartingale.
Applying Ito's formula formally, ignoring the singularity in $|x-y|^{-\rho}$, 
and writing $dM$ for any local martingale terms, we find
\begin{eqnarray}  \label{6.1}
dS_t^{(\rho)} &=& dM+ \int \int u_t(dx)u_t(dy) \left( \frac12 \Delta
\left(|x-y|^{-\rho}\right) + \kappa^2 |x-y|^{-(\rho+2)} \right) dt  \nonumber \\
&=& dM+\left(\rho^2-(d-2)\rho+\kappa^2\right) \int \int
u_t(dx)u_t(dy)|x-y|^{-(\rho+2)} dt  \nonumber
\end{eqnarray}
where $\Delta$ is the Laplacian on $\mathbf{R}^{2d}$, acting on both variables 
$x$ and $y$. The solution to the inequality  $\rho^2- (d-2)\rho+\kappa^2 \leq 0$ gives 
the condition  $\alpha \leq \rho \leq d-2-\alpha$.
The rigorous calculation below does not quite apply to  the boundary value of
$\rho=d-2-\alpha$. 

First we prove a lemma
extending the martingale problem to test functions on $\mathbf{R}^{2d}$.
\begin{Lemma} \label{mgpextension} 
Suppose that $\{u_t(dx)\} $ is a solution to (\ref{1.1}) with initial condition $\mu$.
Then for twice differentiable function $f:  \mathbf{R}^{2d} \rightarrow \mathbf{R}$ with 
compact support 
\begin{equation}
M_t(f) = \int \! \int f(x,y) u_t(dx) u_t(dy)   - \int^t_0 \! \int \! \int \left( \frac12 \Delta f(x,y) 
+ \kappa^2  \frac{f(x,y)}{|x-y|^2} \right) u_s(dx) u_s(dy) ds
\end{equation}
defines a continuous local martingale.
\end{Lemma}

\noindent
\textbf{Proof.} For $f$ of product form, that is $f(x,y) = \sum_{k=1}^n \phi_k(x) \psi_k(y)$ 
where $\phi_k,\psi_k \in C^2_c$, this claim is a consequence of the martingale problem 
(\ref{1.4}) and (\ref{1.5})  together with integration by parts. Now we claim that we can choose $f_n(x,y)$ 
of product form, and with a common compact support, so that $f_n$ and $\Delta f_n$ converge 
uniformly to $f$ and $\Delta f$. One way to see this is consider the one point compactification 
$E$ of the open box  $\{(x,y): |x|,|y| < N\}$ and to let  $(X_t,Y_t)$ be independent d-dimensional 
Brownian motions absorbed on hitting the boundary point of $E$.  Then consider  the algebra 
$\mathcal{A}$  generated by the constant functions and the product functions $\phi(x) \psi(y)$, where 
$\phi,\psi$ are compactly supported in $\{x: |x|<N\}$. The 
Stone-Weierstrass theorem shows that this algebra is dense in the 
space of continuous functions on $E$ and the transition semigroup $\{T_t\}$ of $(X_t,Y_t)$ maps 
$\mathcal{A} $ to itself. A lemma of Watanabe (see \cite{Ethier+Kurtz86} Proposition 3.3) now 
implies that $\mathcal{A}$ is a core for the generator of $(X_t,Y_t)$  and this implies the above claim.

The continuity of $t \to u_t$, and the calculation in
Lemma \ref{momentcor1} part ii), imply that  $M_t(f_n)$ converges to $M_t(f)$ 
uniformly on compacts, in probability.  So the limit $M_t(f)$ has continuous
paths. Also, if $f_n$ and $f$ are supported in the compact set $A$,  the stopping times 
\[
T_k = \inf \{t: u_t(A) + \int^t_0 \int_{A} \int_{A} 
\frac{ u_s(dx) u_s(dy)}{|x-y|^2} ds \geq k \} 
\]
satisfy $T_k \uparrow \infty$ and reduce all the local martingales $M_t(f_n)$
to bounded martingales. We may then pass to the 
limit as $n \rightarrow \infty$ to see that $M_t(f)$ is a local martingale reduced by $\{T_k\}$.  
\qed

\noindent
\textbf{Proof of Proposition \ref{dimension}.} For part i) we may, by conditioning on the 
initial condition, suppose that 
$u_0 = \mu   \in \mathcal{H}_{(d-2-\alpha)-}$. We may then choose $a$ so that
$\mu \in  \mathcal{H}^a_{\beta}$ for all $\beta < d-2-\alpha$. 

We shall approximate $|x-y|^{-\rho}$ by a sequence of compactly
supported functions as follows. Choose $\phi_n \in \mathcal{C}_c ^2$
satisfying  and $\phi_n (x) = 1$ for $|x| \leq n$ and with $\phi_n, \partial_{x_i} \phi_n,
\partial_{x_i x_j} \phi_n$ uniformly bounded over $x$ and $n$. 
Define, for $\varepsilon>0$  and $\rho \in [(d-2)/2, d-2-\alpha]$, 
\[
f_{n, \varepsilon}(x,y) = \left( 1+ ( \varepsilon + |x-y|^2 )^{-\rho/2} \right)
e^{ - a (1+|x|^2)^{1/2} - a(1+|y|^2)^{1/2}} \phi_n(x) \phi_n (y). 
\]
A calculation shows that 
\begin{eqnarray*}
\lefteqn{\left( \frac12 \Delta + \frac{\kappa^2}{|x-y|^2} \right) (\varepsilon +
|x-y|^2)^{-\rho/2}} \\
& = & \left( \varepsilon + |x-y|^2 \right)^{-(\rho+4)/2} 
\left( (\rho^2-(d-2)\rho+\kappa^2) |x-y|^2 + ( 2\kappa^2
-\rho d) \varepsilon + \varepsilon^2 \kappa^2 |x-y|^{-2} \right) \\
& \leq & \left( \varepsilon + |x-y|^2 \right)^{-(\rho+4)/2} \varepsilon^2
\kappa^2 |x-y|^{-2} \\
& \leq & C(\rho) |x-y|^{-(\rho+2)}.
\end{eqnarray*}
The penultimate inequality follows from the restriction on the value of $\rho$
and $\kappa$ and the last inequality follows by considering separately the
cases $\varepsilon<|x-y|^2$ and $\varepsilon \geq |x-y|^2$. 
A simple calculation also shows that 
\[
\left| \Delta \left(e^{ - a (1+|x|^2)^{1/2} - a(1+|y|^2)^{1/2}}\right) \right| \leq 
C(a) e^{ - a (1+|x|^2)^{1/2} - a(1+|y|^2)^{1/2}}
\leq C(a) e^{-a|x|-a|y|}.
\]
Now a lengthy calculation, using the above two bounds as key steps, shows that
\[
\left| \left( \frac12 \Delta + \frac{\kappa^2}{|x-y|^2} \right) f_{n, \varepsilon}
(x,y) \right| \leq C(a,\rho) \left( 1 + |x-y|^{-(\rho+2)}\right) e^{-a|x|-a|y|}. 
\]
Note the bound is uniform over $n$ and $\epsilon$.
Using the test function $f_{n, \varepsilon}(x,y)$ in Lemma \ref{mgpextension}
we have that 
\begin{equation}  \label{6.2}
M_t(f_{n, \varepsilon}) = \int \int f_{n,\varepsilon}(x,y) u_t(dx) u_t(dy) - 
\mathcal{E}(\varepsilon,t)
\end{equation}
is a continuous local martingale and 
\begin{equation}  \label{6.3}
\mathcal{E}(\varepsilon,t) \leq C(a,\rho) \int^t_0 \int_{\mathbf{R}^d} 
\int_{\mathbf{R}^d} \left( 1+| x-y|^{-(\rho+2)}\right) e^{-a|x| -a|y|} u_s(dx) u_s(dy) ds.
\end{equation}
Now we apply Doob's inequality in the following form

\begin{Lemma} \label{dooblemma}
Suppose $\{A_t\}, \, \{M_t\}, \, \{D_t\}$ are continuous processes satisfying $0 \leq A_t = M_t +D_t$ 
and where $M_t$ is a continuous local martingale with $M_0$ bounded. Then for $\lambda \geq 0$ 
\[
P \left( \sup_{s \leq t} A_s \geq 2 \lambda \right) \leq \frac{1}{\lambda}
\left( E[M_0]+3E [ \sup_{s \leq t}|D_s| ] \right). 
\]
\end{Lemma}

\noindent
\textbf{Proof} If $\{T_k\}$ reduce the local martingale $M_t$ then by Doob's
inequality for positive submartingales 
\begin{eqnarray*}
P(\sup_{s \leq t \wedge T_k}|M_s| \geq \lambda) 
& \leq & \frac{1}{\lambda} E[|M_{t \wedge T_k}|] \\
& \leq & \frac{1}{\lambda} \left( E[A_{t \wedge T_K}] + E [|D_{t \wedge T_K}|] \right) \\
& \leq & \frac{1}{\lambda} \left( E[M_0] + 2 E[|D_{t \wedge T_K}|] \right) \\
& \leq & \frac{1}{\lambda} \left( E[M_0] + 2 E[ \sup_{s \leq t} |D_s|] \right) .
\end{eqnarray*}
Let $k \rightarrow \infty$ and combine with the bound 
$P( \sup_{s \leq t} |D_s| \geq \lambda) \leq E [ \sup_{s \leq t} |D_s| ] / \lambda$ to complete
the lemma.  \qed

We apply this lemma to the decomposition (\ref{6.2}) together with the
bound (\ref{6.3}) to obtain 
\begin{eqnarray}   
\lefteqn{ P \left( \sup_{s \leq t} \int \int 
 \left( 1+ |x-y|^{-\rho} \right) e^{-a|x|-a|y|}  u_t(dx) u_t(dy)  > 2 \lambda \right) } 
\nonumber \\
& = & \lim_{\varepsilon \rightarrow 0, n \to \infty} P \left( \sup_{s \leq t} \int \int
f_{n,\varepsilon} (x,y) u_t(dx) u_t(dy) > 2 \lambda \right)  \nonumber \\
& \leq & \frac{C(a,\rho)}{\lambda}  \|\mu(dx) \exp(-a|x|)  \|_{\rho}^2 \nonumber \\
&& \hspace{.2in} + \frac{C(a,\rho)}{\lambda}   E \left[ \int^t_0 \int \int
\left( 1+  |x-y|^{-(\rho+2)}\right)  e^{-a|x|-a|y|} u_s(dx) u_s(dy) ds \right]. \label{6.4}
\end{eqnarray}
A little effort, as in Lemma  \ref{momentcor1} part ii) and using the fact that
$\mu \in  \mathcal{H}^a_{\beta}$ for all $\beta < d-2-\alpha$, 
shows that the expectation on the right hand side of
(\ref{6.4})  is finite.  One needs, however, the strict inequality  $\rho < d-2-\alpha$ 
so that the worst pole  is $|x'-y'|^{-(\rho+2+\alpha)}$, which is therefore still integrable 
ensuring the bound (\ref{1.12}) applies. The bound in (\ref{6.4})
implies part i) of the Proposition.

For part ii) we may suppose, by conditioning 
on the initial condition, that $u_0 = \mu \in \mathcal{H}^0_{\alpha+}$. 
But then Lemma \ref{momentcor1} part ii) implies, for fixed $t_0 >0$, that 
$u_{t_0}(dx) \in \mathcal{H}_{d-2-\alpha}$ almost surely. The Markov 
property of solutions and part i) then imply that the desired conclusion holds for 
$t \geq t_0$. Letting  $t_0 \downarrow 0$ completes this proof. \qed

\begin{Corollary} \label{strongmarkovproperty}
The family $ \{ Q_{\mu}: \mu \in \mathcal{H}_{(d-2-\alpha)-}\}$ is a strong Markov family.
\end{Corollary}

\noindent
\textbf{Proof.}
Let $\{u_t\}$ be a solution defined on $(\Omega, \mathcal{F}, \mathcal{F}_t, P)$ and
satisfying $ P(u_0 \in \mathcal{H}_{(d-2-\alpha)-})=1$. 
Let $\tau<\infty $ be a $ \mathcal{F}_t$ stopping time and $\tau(n) < \infty $ be
discrete stopping times satisfying $\tau_n \downarrow \tau$. Fix
$0 \leq t_1 < \ldots < t_n$ and $f_1, \ldots , f_n \in \mathcal{C}_c$ and set
$H(u) = \exp(i ( u_{t_1}(f_1) + \ldots + u_{t_n}(f_n)))$. Fix a set 
$\Lambda \in \mathcal{F}_{\tau}$.   Then the ordinary Markov property 
implies that 
\begin{equation} \label{6.5}
E\left[ H( u_{\tau(n) + \cdot}) \mathbf{1}(\Lambda) \right] 
= E \left[ Q_{u_{\tau(n)}}[H] \mathbf{1} (\Lambda) \right].
\end{equation}
If we can pass to the limit as $n \to \infty$ to replace 
$\tau (n)$ by $\tau$, then this identity implies the result.  By the continuity
of paths the left hand side of (\ref{6.5}) converges as desired. We claim that 
\[
\mbox{if $\mu_n \to \mu$ vaguely and $\sup_n \| \mu_n (dx) e^{-a|x|} \|_{\alpha} < \infty$
then $Q_{\mu_n}[H] \to Q_{\mu}[H]$.}
\]
Assuming  this claim, Proposition \ref{dimension} part i) 
allows us to pass to the limit on the right hand side of (\ref{6.5}).
To prove the claim we let $u_t(dx)$ be the solution starting at 
$\mu$ constructed using the chaos expansion and 
$u_{N,t}$ the approximation using only the first $N$ terms of the expansion.
Then 
\begin{eqnarray*}
Q_{\mu}[H] & =&  E\left[ \exp( i \sum_{j=1} u_{t_j}(f_j)) \right] \\
& =&  E\left[ \exp( i \sum_{j=1} u_{N,t_j}(f_j)) \right] 
+ \mbox{Error}(N,\mu)
\end{eqnarray*}
where 
\[
\left| \mbox{Error}(N,\mu) \right| \leq \left( \sum_{j=1}^n E\left [(u_{t_j}(f_j)-u_{N,t_j}(f_j))^2 \right] \right)^{1/2}.
\]
The function $ E\left[ \exp( i \sum_{j=1} u_{N,t_j}(f_j)) \right]$ is continuous in
$\mu$ and $\mbox{Error}(N,\mu) \to 0$ as $N \to \infty$.
So the claim follows if we can show $\sup_n |\mbox{Error}(N,\mu_n)| \to 0$ as $ N \to \infty$.
Using the isometry as in Lemma \ref{series-converges} we see that
\[
E\left [(u_{t}(f)-u_{N,t}(f))^2 \right] = \int \int H_N(x,y) \mu(dx) \mu(dy)
\]
where 
\[
H_N(x,y) = \int \int f(x') f(y') G_t(x-x') G_t (y-y') E_{0,x,y}^{t,x',y'}
\left[ \sum_{k=N+1}^{\infty} \frac{1}{k!} \left( \int^t_0 \frac{\kappa^2 ds}{|X^1_s-X^2_s|^2} \right)^k
\right] dx' dy'
\]
is bounded by 
\[
H_N(x,y) \leq C(a, t) e^{-a|x|-a|y|} (1+ |x-y|^{-\alpha}).
\]
Note that $H_N(x,y)$ is monotone decreasing but not continuous. 
The assumptions of the claim allow, by an approximation argument, to ignore the 
singularity in the function $H_N(x,y)$ and replace it  by a 
monotone decreasing continuous function $\tilde{H}_N(x,y)$ of compact support.
But then the vague convergence $\mu_n \to \mu$ implies that 
$\sup_n \int \int \tilde{H}_N(x,y) \mu_n(dx) \mu_n(dy) \downarrow 0$ as $N \to \infty$
(for example by the argument of Dini's lemma). This completes the proof of the claim.
\qed
%
%
\subsection{Density of Support} \label{s7.2}
%
%
In this subsection we give the proof of Theorem \ref{t5} ii).  We start with
an outline of the method. Assume that 
$u_0 ( B_r(a) )>0 $ and fix $T>0$. We wish to show that with probability one
$u_{T} ( B_r(b) )>0 $. We consider various tubes in $[0,T] \times \mathbf{R}^d$
which connect $\{0\}\times B_r(a)$ with $\{T\}\times B_r(b)$. (By a tube we
mean that for any time $t$ the cross section of the tube with the slice 
$\{t\}\times \mathbf{R}^d$ is a ball of radius $r$.) We consider a subsolution
to the equation which has Dirichlet boundary conditions on the edge of the tube.
We will show that the probability that the subsolution is non-zero at time $T$
is a constant not depending on the tube. It is possible to construct an infinite family 
of such tubes such that each pair has very little overlap. Then a zero-one law will 
guarantee that, with probability one, at least one of the 
subsolutions will be non-zero. Applying this for a countable family of 
open balls we shall obtain the density of the support. Note this  implies 
that the solution never dies out completely.  Note also that for the equation
(\ref{1.1}) posed on a  finite region the above argument fails, as there
is not enough room to fit an infinite family of nearly disjoint tubes. 

Let us give a rigorous definition of the tubes described above. For a
piecewise smooth function $g: [0,T] \to \mathbf{R}^d$ the tube centered on $g$
is defined as 
\[
\mathbf{T}= \left\{(t,x)\in [0,T] \times \mathbf{R}^d: x \in
B_r(g(t)) \right\}. 
\]
If $\mathbf{T}$ is such a tube, let $\partial \mathbf{T}$ be the boundary of 
$\mathbf{T}$, minus the part of the boundary at $t=0$ and $t=T$. 
Now we aim find a solution $(u^{\mathbf{T}}_t(dx): 0 \leq t \leq T)$
to the equation (\ref{1.1}) but restricted to the tube
$\mathbf{T}$ and with Dirichlet boundary conditions, that is
\begin{equation} \label{6.6}
\left\{ \begin{array}{rcll}
\frac{\partial u^{\mathbf{T}}_t}{\partial t} & = & \Delta  u^{\mathbf{T}}_t+\kappa 
u^{\mathbf{T}}_t\dot F(t,x) & \mbox{for $(x,t) \in \mathbf{T}$,} \\
u^{\mathbf{T}}_0(dx)  & = & \nu (dx), & 
 \mbox{where supp$(\nu) \subseteq B_r(g(0))$,}  \\
u^{\mathbf{T}}_t(dx)  & = & 0 &
\mbox{for $(x,t) \in \partial \mathbf{T}$.} 
\end{array} \right.
\end{equation}
As in Section \ref{s4}, a chaos expansion with respect to the noise $F$ yields 
solutions to (\ref{6.6}). We do not give the proof.
The only changes needed are that the stochastic integrals are restricted to the
tube and the Green's function $G_t(x-y)$ must be replaced by the
Green's function for the tube $G^{\mathbf{T}}_t (x-y)$, that is the fundamental solution for the
heat equation in the tube with Dirichlet boundary conditions.
As in Section \ref{s4} the convergence of the series is guaranteed by 
the finiteness of an exponential Brownian bridge moments; however the moments that 
are needed are of the form
\[ 
E_{0,x,y}^{t,x',y'}\left[\exp\left( \int_{0}^{t}  \frac{\kappa^2}{ \left| X_s^1-X_s^2 \right|^2}
ds \right)  I( \sup_{s \leq t} |X^1_s-g(s)| \vee |X^2_s-g(s)| < r) \right] 
\]
and so are less than the corresponding moments needed to ensure the 
solution on the whole space converges. 

Fix the noise $F$, on its filtered probability space, and construct via chaos expansions
$\{u_t(dx)\}$  the solution to (\ref{1.1}) started at $\mu \in \mathcal{H}_{\alpha+}$ and 
$\{u^{T}_t(dx)\}$  the solution to (\ref{6.6}) started at $\nu= \mu |_{B_r(g(0))}$.
We may also construct approximating solutions $u^{(\mathbf{T},\varepsilon)}_t(x)dx$ to 
$ u^{\mathbf{T}}_t(dx)$, by using the smoother noise $F^{\varepsilon}$ and the initial condition
$\nu^{(\varepsilon)} = G^{\mathbf{T}}_{\varepsilon} \nu$, exactly as we approximated
 $u_t(dx)$ by $u^{(\varepsilon)}_t(x)dx$. A fairly standard comparison argument shows that 
$u^{(\mathbf{T},\varepsilon)}_t(x) \leq  u^{(\varepsilon)}_t(x)$. Passing to the limit as 
$\varepsilon \to 0$ we find that, with probability one, 
\begin{equation} \label{6.7}
u^{\mathbf{T}}_t(dx) \leq  u_t(dx) \qquad \mbox{ for all $ 0 \leq t \leq T$. }
\end{equation}

We now start the proof of Theorem \ref{t5} ii).  As described above it is enough to 
assume that $u_0(B_r(a))>0$, for some $a \in \mathbf{R}^d$ and $r>0$, and 
to show, for fixed $b\in\mathbf{R}^d$, that $u_T(B_r(b))>0$ with probability one. 
For notational ease we shall take $a=b=0$ and $r=1$  since the proof needs only small 
changes  for other values of $r,a,b$. Let $e_1$ be the unit vector $(1,0,\ldots,0)$.  We consider a 
sequence of  piecewise linear functions $g_n(t)$ given, for  $n \in \mathbf{Z}$,  by
\[
g_n(t)= \left\{ 
\begin{array}{ll} 2nte_1 & \mbox{if $\;\; 0\le t \leq T/2$} \\ 
2n(T-t)e_1 & \mbox{if $\;\; T/2 \leq t \leq T$} \end{array} \right. 
\]
We write $\mathbf{T}_n$ for the tube centered on $g_n$. The Feynman-Kac representation 
(\ref{1.18}), adapted for the Dirichlet boundary conditions, gives the following representation for the solution 
$u^{(\mathbf{T}_n, \varepsilon)}_T(f)$ for a test function $f \geq 0 $ supported in $B_1(0)$:
\begin{eqnarray*} 
u^{(\mathbf{T}_n, \varepsilon)}_T(f)
&  = &  e^{-\Gamma_{\varepsilon}(0) T}  \int dx  \int \nu^{(\varepsilon)}(dy)  G_T(x- y) f(x)  \\
&& \hspace{.25in} \cdot  E_{0,y}^{T,x} \left[ \exp \left(  \kappa \int^T_0 F^{(\varepsilon)}(ds,X_{s}) \right)  
I( (s,X_s) \in \mathbf{T}_n, \; \forall s \leq T) \right],
\end{eqnarray*}
where all the integrals are over the ball $B_1(0)$. By conditioning on the position of the Brownian 
bridge at time $T/2$ we find 
\begin{eqnarray} 
 u^{(\mathbf{T}_n, \varepsilon)}_T(f)  
&  = & e^{-\Gamma_{\varepsilon}(0) T} \int dx \int \nu^{(\varepsilon)}(dy) 
\int dz f(x) \label{6.8} \\ 
&& \hspace{.2in} \cdot  G_{T/2}(x- nTe_1 - z) G_{T/2}(nTe_1 + z - y) f(x)  E_1(y,z) E_2(z,x)
\nonumber
\end{eqnarray}
where 
\[
E_1(y,z) =  E_{0,y}^{T/2,z+nT} \left[ \exp \left(  \kappa \int^{T/2}_0 F^{(\varepsilon)} (ds,X_{s}) \right)  
I( (s,X_s) \in \mathbf{T}_n, \;  \forall s \leq T/2) \right] 
\]
and 
\[ 
E_2(z,x) =  E_{0,z+nT}^{T/2,x} \left[ \exp \left(  \kappa \int^{T}_{T/2} F^{(\varepsilon)} (ds, X_{s-(T/2)}) 
\right)   I( (s,X_{s-(T/2)}) \in \mathbf{T}_n, \; \forall s \leq T) \right] 
\]
All the randomness in the representation (\ref{6.8}) is contained in the Brownian bridges
expectations $E_1(x,z) $ and $E_2(z,y)$.  By adding a suitable linear drift to 
the Brownian bridge we may rewrite 
\[ 
E_1(y,z) =  E_{0,y}^{T/2,z} \left[ \exp \left(  \kappa \int^{T/2}_0 F^{(n, \varepsilon)} (ds,X_{s}) \right)  
I( (s,X_s) \in \mathbf{T}_n, \; \forall s \leq T/2) \right] 
\]
where $F^{(n, \varepsilon)}(x,t) = F^{(\varepsilon)} (x+nt,t)$ is a new noise which has the same
covariance structure as $F^{(\varepsilon)}$. This shows that the laws of $E_1(y,z)$ is independent 
of $n$, and a similar argument applies to $E_2(z,x)$ which is also independent of $E_1(y,z)$. 
 Also  for $x,z \in B_1(0)$ 
\[ 
\frac{G_{T/2} (x- nTe_1-z)}{G_{T/2}(x-z)} = \exp(-n^2T - 2ne_1 \cdot (x-z)) \geq \exp(-n^2T - 4|n|). 
\]
A similar lower bound holds for $G_{T/2}(nTe_1+z-y)$. Using these in  (\ref{6.8}) we see that
the variable $ u^{(\mathbf{T}_n,\varepsilon)}_T(f) $ stochastically dominates the variable 
$ C(n,T)   u^{(\mathbf{T}_0,\varepsilon)}_T(f)$, 
where $C(n,T)$ is a strictly  positive constant independent of $\varepsilon$. Letting $\varepsilon 
\downarrow 0$ we obtain  the same  stochastic dominance for the solutions driven by the 
singular noise $F$:
\[
u^{(\mathbf{T}_n)}_T(f) \stackrel{s}{\geq} C(n,T)  u^{(\mathbf{T}_0)}_T(f),
\]
where the inequality stands for stochastic domination. Let $A_n$ be the event  
$ \{ u_T^{\mathbf{T}_n} (B_1(0))>0 \}$.   Then $P(A_n) \geq P(A_0)$ 
for all $n$ by this stochastic domination. Also $P(A_0)>0$, as can be seen from the fact that the first 
moment of $  u_T^{(n)} (B_1(0))$ is given, in a similar way as for the first moments
in (\ref{1.6}), by the heat flow in the tube and is hence non-zero.

Finally we apply a zero-one law to conclude the result. Consider the sequence of noises defined by
\[
F_n = \left( \dot{F}(t,x +g_n(t)) : 0<t<T, \, |x| <1 \right) \quad \mbox{for n=0,1,\ldots}
\]
Since the correlation structure of $F$ is unchanged by piecewise linear
shifts the noises $\{F_k\}$ are identically distributed and form a stationary
sequence. We claim this sequence is also strong mixing. For this it is enough to show,
for all $k$ and bounded measurable $G,H$, that as $n \in \mathbf{Z}$
\begin{equation} \label{6.9}
E\left[ G(F_{-k}, \ldots, F_k) H(F_{n-k}, \ldots, F_{n+k}) \right]
\to E\left[ G(F_{-k}, \ldots, F_k)\right]  E \left[ H(F_{-k},  \ldots, F_{k}) \right] . 
\end{equation}
Suppose that $\phi_{i,j}(x,t)$ are test functions supported
in $(0,T) \times B_1(0)$. Suppose $G$ and $H$ are bounded continuous functions
of the vector 
\[ 
\left(\int_0^T \int_{B_1(0)} \phi_{i,j} dF_i: -k \leq i \leq k, j=1,
\ldots, k \right). 
\]
Each integral $\int \int \phi_{i,j} dF_i$ is a Gaussian variable. Also, the covariance 
between $\int \int \phi_{i,j} dF_i$ and $\int \int \phi_{i+n,j'} dF_{i+n}$
converges to zero as $n \to \infty$. This implies that the mixing relation
(\ref{6.9}) holds for $G,H$ of this special. A monotone class argument then
proves the mixing relation for general $G$ and $H$. 

Define $\mathcal{S}_n$ to be the $\sigma$-field generated by the noises 
$(F_n, F_{n+1}, F_{n+2}, \ldots)$.  The strong mixing of the sequence
implies that the sigma field 
$ \mathcal{S} = \cap_{n=1}^{\infty} \mathcal{S}_n $ is trivial in that 
$P(S)=0$ or $1$ for all $S \in \mathcal{S}$. 
The construction of the solutions by a Wiener chaos expansion shows that
the solution $u^{(\mathbf{T}_n)}$ is measurable with respect to
the sigma field generated by the noise $ \dot{F}(t,x)$ for $(t,x) \in T_n$.
Thus  the event  $A_n  $ is $\mathcal{S}_n$ measurable and the event 
$\{A_n \; i.o.\}$ is $\mathcal{S}$ measurable. Since $P(A_n) $ is bounded below 
uniformly in $n$ the event $\{A_n \; i.o.\}$ must have probability one. Finally, since
 $u_T(B_1(0)) \geq \sup_n u^{(\mathbf{T}_n)}_T(B_1(0))$ by (\ref{6.7}) the
proof is complete.  
%
%
\subsection{Singularity of solutions} \label{s7.3}
%
%
In this subsection we prove the singularity assertion in Theorem \ref{t5} iii). We first 
sketch a short argument that suggests the solutions are singular. Fix $T>0$ and 
$x\in\mathbf{R}^d$. For $t\in [0,T)$ we consider the process 
\[
M_t(x)=\int G_{T-t}(x-y) u_t(dy). 
\]
It is possible to extend  the martingale problem (\ref{1.4}) to test functions
that depend on time and that do not have compact support,  provided that they decay faster than 
exponentially at infinity. Using the test function  $(t,y) \rightarrow G_{T-t}(x-y)$ it follows from this 
extension that $\{M_t\}$ is a nonnegative  continuous local martingale for $t \in [0,T)$. The explosion 
principle (see \cite{Rogers+Williams00} Corollary IV. 34.13) implies that the quadratic variation must 
remain bounded as $t \uparrow T$. Therefore, with probability 1, 
\begin{equation}  \label{6.10}
\langle M(x) \rangle_T = \int_{0}^{T}\int \int u_t(dy)u_t(dz)G_{T-t}(x-y)G_{T-t}(x-z)|y-z|^{-2} < \infty.
\end{equation}
However, a short calculation shows that if $u_t(y)$ has a continuous, strictly positive density in the 
neighborhood of $(T,x)$ then the integral in (\ref{6.10}) is infinite.

Instead of pursuing this argument we show that the scaling relation can be used to convert 
the death of solutions at large times to the singularity of solutions at a fixed time. Applying 
the scaling Lemma \ref{scaling}, with the choices $a = \varepsilon^{-d}$, $b=\varepsilon^2$ 
and $c = \varepsilon$, we find that, under the initial condition $u_0(dx) = CL(dx)$ (where $L(dx)$ 
is Lebesgue measure), that $u_t(B(0,\varepsilon)) $ has the same distribution as 
$\varepsilon^{d} u_{t/\varepsilon^2} (B(0,1))$. Also, as in the proof of Theorem \ref{t4} ii), the 
linearity of the equation and the concavity of the function $z \rightarrow \sqrt{z} $ imply that 
the map $\mu \rightarrow Q_{\mu}[U_t(B(0,\varepsilon))^{1/2}]$ is increasing and concave 
in $\mu$.

Take a solution $\{u_t(dx)\}$ with $u_0$ of locally bounded intensity. Then, for fixed $t>0$, 
\begin{eqnarray*}
\lefteqn{E\left[ \frac{u_t(B(0,\varepsilon))^{1/2}}{\varepsilon^{d/2}} \right]} \\
& = & \int_{\mathcal{H}_{\alpha+}} Q_{\mu} \left[ 
\frac{U_{t/2}(B(0,\varepsilon))^{1/2}}{\varepsilon^{d/2}} \right] P(u_{t/2} \in d\mu) 
\quad \mbox{(by the Markov property)} \\
& \leq & Q_{ E(u_{t/2}(dx))} \left[ 
\frac{U_{t/2} (B(0,\varepsilon))^{1/2}}{\varepsilon^{d/2}} \right] \quad 
\mbox{(by Jensen's inequality)} \\
& \leq & Q_{ C(t) L(dx) } \left[ 
\frac{U_{t/2} (B(0,\varepsilon))^{1/2}}{\varepsilon^{d/2}} \right] \quad 
\mbox{(since $u_0$ has bounded intensity)} \\
& = & Q_{ C(t) L(dx) } \left[ U_{t/2 \varepsilon^2} (B(0,1))^{1/2} \right]
\quad \mbox{(by scaling)} \\
& = & \left(C(t)\right)^{1/2} Q_{ I_{B(0,1)} } \left[ U_{t/2 \varepsilon^2} (1)^{1/2}
\right] \quad \mbox{(by self-duality, Proposition \ref{dualityrelation})} \\
& \rightarrow & 0 \quad \mbox{as $\varepsilon 
\downarrow 0 \quad $ (by Proposition \ref{liggett}).}
\end{eqnarray*}
The same result holds true if $B(0,\varepsilon)$ is replaced by $B(x, \varepsilon)$ for 
any $x \in \mathbf{R}^d$. We may decompose the measure $u_t = u_t^{(ac)} + u_t^{(s)}$ 
into its absolutely continuous and singular parts and write $u_t^{(ac)} = A_t(x) dx$ for a 
locally $L^1$ function $A_t(x) \geq 0 $. Then
\begin{eqnarray*}
E\left[ \int A_t^{1/2}(x) dx \right] & = &  \int E \left[ \lim_{\varepsilon \downarrow 0} 
\frac{u^{(ac)}_t(B(x,\varepsilon)^{1/2}}{\varepsilon^{d/2}} \right] dx 
\quad \mbox{(Lebesgue differentiation theorem)} \\
& \leq & \int \lim_{\varepsilon \downarrow 0} E \left[ 
\frac{u^{(ac)}_t(B(x,\varepsilon)^{1/2}}{\varepsilon^{d/2}} \right] dx 
\quad \mbox{(Fatou's lemma)} \\
& \leq & \int \lim_{\varepsilon \downarrow 0} E \left[ 
\frac{u_t(B(x,\varepsilon)^{1/2}}{\varepsilon^{d/2}} \right] dx =0.
\end{eqnarray*}
Thus $A_t = 0$ with probability one.

In general we may decompose the initial condition $u_0 \in \mathcal{H}_{\alpha+}$ 
as a countable sum of measures $u_0 = \sum_n u^{(n)}_0$ where each $u^{(n)}_0$
has locally bounded intensity. Use a single  noise
to define chaos expansion solutions $u^{(n)}_t(dx)$ with initial conditions
$u^{(n)}_0$. It is easy to check that $\sum_n u^{(n)}_t(dx)$ is a solutions started at $u_0$.
Then, applying the above 
argument to each $u_t^{(n)}$ yields the desired result in the general case. This completes 
the proof of Theorem \ref{t5} iii).

\noindent 
\textbf{Remark:} We would like to thank M. Yor for informing us of his work in \cite{Yor80}. 
%
%
\bibliographystyle{plain}
\bibliography{bibtex}

\end{document}